\documentclass[12pt]{article}
\usepackage{amssymb,a4}
\usepackage{amsmath,amsfonts,amssymb,amsthm}
\def\Res{{\rm Res}}
\def\wt{{\rm wt}}

\def\C{{\mathbb C}}

\def\Z{{\mathbb Z}}

\def\1{{\bf 1}}
\def \pf{\noindent {\bf Proof: \,}}
\def\theequation{5.\arabic{equation}}

\renewcommand{\theequation}{\thesection.\arabic{equation}}

\newtheorem{theorem}{Theorem}[section]
\newtheorem{prop}[theorem]{Proposition}
\newtheorem{lem}[theorem]{Lemma}
\newtheorem{coro}[theorem]{Corollary}
\newtheorem{remark}[theorem]{Remark}
\theoremstyle{definition}

\begin{document}
\begin{center}
{\Large {\bf Bimodules associated to  vertex operator algebras}} \\

\vspace{0.5cm} Chongying Dong\footnote{Supported by NSF grants,
China NSF grant 10328102 and  a Faculty research grant from  the
University of California at Santa Cruz.}
\\
Department of Mathematics\\ University of
California\\ Santa Cruz, CA 95064 \\
Cuipo Jiang\footnote{Supported in part by China NSF grant
10571119.}\\
 Department of Mathematics\\ Shanghai Jiaotong University\\
Shanghai 200030 China\\
\end{center}
\hspace{1cm}
\begin{abstract} Let $V$ be a vertex operator algebra and $m,n\geq 0.$
We construct an $A_n(V)$-$A_m(V)$-bimodule $A_{n,m}(V)$ which
determines the action of $V$ from the level $m$ subspace to level $n$
subspace of an admissible $V$-module. We show how to use $A_{n,m}(V)$ to construct naturally
admissible $V$-modules from $A_m(V)$-modules. We also determine the structure
of $A_{n,m}(V)$ when $V$ is rational.   2000MSC:17B69
\end{abstract}

\section{Introduction}
\def\theequation{1. \arabic{equation}}
\setcounter{equation}{0}

The representation theory for a vertex operator algebra \cite{B},
\cite{FLM} has been
studied largely in terms of representation theory for various
associative algebras associated to the vertex operator algebra (see
\cite{Z}, \cite{KW}, \cite{DLM2}-\cite{DLM4}, \cite{MT},
\cite{DZ}, \cite{X}). A sequence of associative algebras $A_n(V)$
for $n\geq 0$ was introduced in \cite{DLM3} to deal with  the
first $n+1$ homogeneous subspaces of an admissible module. These
algebras extend and generalize the associative algebra $A(V)$
constructed in \cite{Z}. The main idea of $A_n(V)$ theory is how
to use the first few homogeneous subspaces of a module to
determine the whole module. From this point of view, the $A_n(V)$
theory is an analogue of the highest weight module theory for
semisimple Lie algebras in the field of vertex operator algebra.

Let $M=\bigoplus_{n=0}^{\infty}M(n)$ be an admissible $V$-module with
$M(0)\ne 0.$ (cf. \cite{DLM2}).Then each $M(k)$ is an $A_n(V)$-module
for $k\leq n$ \cite{DLM3}.
 On the other hand, given an $A_n(V)$-module
$U$ which cannot factor through $A_{n-1}(V)$ one can construct a Verma
type admissible $V$-module $\bar M(U)$ such that $\bar M(U)(n)=U.$ Also $V$ is rational if and only if $A_n(V)$ is
 semisimple for all $n.$ So the collection of associative algebras $A_n(V)$ determine the
 representation theory of $V$ in some sense.
However, $A_n(V)$ preserves each homogeneous subspace
$M(m)$ for $m\leq n$ and cannot map  $M(s)$ to $M(t)$ if $s\ne t.$
The goals of the present paper are to alleviate this situation.

Given two nonnegative integers $m,n$ we will construct an
$A_n(V)$-$A_m(V)$-bimoulde $A_{n,m}(V)$ with the property that for
any $A_m(V)$-module $U$ which cannot factor through $A_{m-1}(V)$
one can associate a Verma type admissible $V$-module
$M(U)=\bigoplus_{k=0}^{\infty}M(U)(n)$ such that
$M(U)(n)=A_{n,m}(V)\otimes _{A_m(V)}U.$ The action of $V$ on
$M(U)(n)$ is determined by a canonical bimodule homomorphism from
$A_{p,n}(V)\otimes_{A_n(V)}A_{n,m}(V)$ to $A_{p,m}(V).$ Also, for
a given admissible $V$-module $W=\bigoplus_{k\geq 0}W(k)$ with
$W(0)\ne 0,$ there is an $A_n(V)$-$A_m(V)$-bimodule homomorphism
from $A_{n,m}(V)$ to ${\rm Hom}_{\C}(W(m),W(n)).$ So the
collection of $A_{n,m}(V)$ for all $m,n\in\Z_+$ determine the
action of vertex operator algebra $V$ on its admissible module $W$
completely. This, in fact, is our original motivation to define
$A_{n,m}(V).$

If $V$ is a rational vertex operator algebra, then $V$ has only
finitely many irreducible admissible $V$-modules up to isomorphism
and each irreducible $V$-module is ordinary (see \cite{DLM2}). In
this case we let $W^1,...,W^s$ be the inequivalent irreducible
admissible $V$-modules such that $W^i(0)\ne 0.$ Then $A_n(V)$ is
the direct sum of full matrix algebras
$A_n(V)=\bigoplus_{i=1}^s\bigoplus_{k=0}^n{\rm End}_{\C}(W^i(k)).$
We show in this paper that if $V$ is rational then
$$A_{n,m}(V)\cong \bigoplus_{l=0}^{{\rm
min}\{m,n\}}\left(\bigoplus_{i=1}^{s}{\rm Hom}_{\mathbb
C}(W^i(m-l),W^i(n-l))\right).
$$
The structure of $A_{n,m}(V)$ for general $V$ will be studied in a sequel
to this paper.

We have already mentioned the $A_n(V)$ theory. In fact, the Verma
type admissible $V$-module $M(U)$ has been constructed
and denoted by $\bar M(U)$ in \cite{DLM3}
using the idea of induced module in Lie theory. But
our work in this paper leads to a strengthening of this old
construction. While the
old construction in \cite{DLM3} was given as an abstract quotient
of certain induced module for certain Lie algebras, the new
construction is explicit and each homogeneous subspace $M(U)(n)$
is obvious. In the case that $U=A_m(V)$ we see immediately that
$\bigoplus_{n\geq 0}A_{n,m}(V)$ is an admissible $V$-module. We
expect that the study of bimodules $A_{n,m}(V)$ will lead to a
proof of some well known conjecture in representation theory.

The result in this paper is also related to some results obtained
in \cite{MNT} where the universal enveloping algebra of a vertex
operator algebra is used instead of vertex operator algebra itself
in this paper.

The paper is organized as follows: In Section 2 we introduce the
$A_n(V)$-$A_m(V)$-bimodule $A_{n,m}(V)$ with lot of technical
calculations. In Section 3 we discuss the various properties of
$A_{n,m}(V)$ such as the $A_{n}(V)$-$A_m(V)$-bimodule epimorphism
from $A_{n,m}(V)$ to $A_{n-1,m-1}(V)$ induced from the identity
map on $V,$ isomorphism between $A_{n,m}(V)$ and $A_{m,n}(V)$ and
bimodule homomorphism from $A_{n,p}(V)\otimes_{A_p(V)}A_{p,m}(V)$
to $A_{n,m}(V).$ In Section 4 we first give an
$A_n(V)$-$A_m(V)$-bimodule homomorphism from $A_{n,m}(V)$ to ${\rm
Hom}_{\C}(W(m),W(n))$ for any admissible $V$-module
$W=\bigoplus_{k=0}^{\infty}W(k).$ We also show how to construct an
admissible $V$-module from an $A_m(V)$-module which cannot factor
through $A_{m-1}(V)$ by using $A_{n,m}(V).$ In addition we show
that $A_n(V)$ and $A_{n,n}(V)$ are the same although $A_{n,n}(V)$
as a quotient space of $A_{n}(V)$ seems much smaller from the
definition. The explicit structure of $A_{n,m}(V)$ is determined
if $V$ is rational.

We assume that the reader is familiar with the basic knowledge on
the representation theory such as of weak modules,
admissible modules and (ordinary) modules as presented in  \cite{DLM1}-\cite{DLM2}
(also see \cite{FLM}, \cite{LL}).

\section{$A_{n}(V)$-$A_{m}(V)$-bimodule $A_{n,m}(V)$}
\def\theequation{2.\arabic{equation}}
\setcounter{equation}{0}

Let $V=(V,Y,1,\omega)$ be a vertex operator algebra. An
associative algebra $A_n(V)$ for any nonnegative integer $n$ has
been constructed in \cite{DLM3} to study the representation theory
for vertex operator algebras (see below). For $m,n\in{\mathbb
Z}_{+}$, we will construct an $A_{n}(V)$-$A_{m}(V)$-bimodule
$A_{n,m}(V)$ in this section. It is hard to see in this section why
$A_{n,m}(V)$ is so defined and the motivation for defining
$A_{n,m}(V)$ comes from the representation theory of $V$
(see Section 4 below).

For homogeneous $u\in V,$ $v\in V$ and $m,n,p\in{\mathbb
Z}_{+}$, define the product $\ast_{m,p}^{n}$  on $V$ as follows
$$
u\ast_{m,p}^{n}v=\sum\limits_{i=0}^{p}(-1)^{i}{m+n-p+i\choose
i}{\rm Res}_{z}\frac{(1+z)^{wtu+m}}{z^{m+n-p+i+1}}Y(u,z)v.
$$
If $n=p$, we denote $\ast_{m,p}^{n}$ by $\bar{\ast}_{m}^{n}$, and
if $m=p$, we denote $\ast_{m,p}^{n}$ by $\ast_{m}^{n}$, i.e.,
$$
u\ast_{m}^{n}v=\sum\limits_{i=0}^{m}(-1)^{i}{ n+i\choose i}{\rm
Res}_{z}\frac{(1+z)^{wtu+m}}{z^{n+i+1}}Y(u,z)v,
$$$$
u\bar{\ast}_{m}^{n}v=\sum\limits_{i=0}^{n}(-1)^{i}{ m+i\choose
i}{\rm Res}_{z}\frac{(1+z)^{wtu+m}}{z^{m+i+1}}Y(u,z)v.$$
The products $u\ast_{m}^{n}v$ and $
u\bar{\ast}_{m}^{n}v$ will induce the right $A_m(V)$-module
and left $A_n(V)$-module structure on $A_{n,m}(V)$ which will
be defined later.

If $m=n,$
then  $u\ast_{m}^{n}v$ and $u\bar{\ast}_{m}^{n}v$ are
equal, and have been defined in \cite{DLM3}. As in \cite{DLM3} we
will denote the product by $u\ast_{n}v$ in this case.

Let $O'_{n,m}(V)$ be the linear span of all $u\circ_{m}^{n}v$ and
 $L(-1)u+(L(0)+m-n)u$, where for homogeneous $u\in V$ and $v\in V$,
$$
u\circ_{m}^{n}v={\rm
Res}_{z}\frac{(1+z)^{wtu+m}}{z^{n+m+2}}Y(u,z)v.$$ Again if $m=n,$
$u\circ_{m}^{n}v$ has been defined in \cite{DLM3} where it was
denoted by $u\circ_n v.$ Let $O_{n}(V)=O'_{n,n}(V)$ (see
\cite{DLM3}). The following theorem is obtained in \cite{DLM3}.

\begin{theorem}\label{t2.1} The $A_n(V)$ is an associative algebra with product
$*_n$ with identity $\1+O_n(V).$
\end{theorem}

We will present more results on $A_n(V)$ and its connection with the
representation theory of $V$ from \cite{DLM3} later on when necessary.
In order to define $A_{n,m}(V)$ we need several lemmas. In the case that
$m=n$ most of these lemmas have been proved in \cite{DLM3}. But when $m\ne n$
even if the old proofs given in \cite{DLM3} work but they are much more
complicated and need a lot of modifications. Sometimes we need to find
totally new proofs.

\begin{lem}\label{l2.2} For any $u,v\in V,$ $u\circ_m^nv$ lies in
$O_{n}(V)\bar{\ast}_{m}^{n}V.$

\end{lem}

\pf Let $u,v\in V.$ Then
\begin{eqnarray*}
& &(L(-1)u)\bar\ast_{m}^{n}v=\sum\limits_{i=0}^{n}(-1)^{i}
{m+i\choose i}{\rm
Res}_{z}\frac{(1+z)^{wtu+1+m}}{z^{m+i+1}}\frac{d}{dz}Y(u,z)v\\
& &=-(wtu+m+1)\sum\limits_{i=0}^{n}(-1)^{i}{ m+i\choose i}{\rm
Res}_{z}\frac{(1+z)^{wtu+m}}{z^{m+i+1}}Y(u,z)v\\
& &+\sum\limits_{i=0}^{n}(m+i+1)(-1)^{i}{ m+i\choose i}{\rm
Res}_{z}\frac{(1+z)^{wtu+m+1}}{z^{m+i+2}}Y(u,z)v\\
& &=-(wtu+m+1)\sum\limits_{i=0}^{n}(-1)^{i}{ m+i\choose i}{\rm
Res}_{z}\frac{(1+z)^{wtu+m}}{z^{m+i+1}}Y(u,z)v\\
& &+\sum\limits_{i=0}^{n}(m+i+1)(-1)^{i}{ m+i\choose i}{\rm
Res}_{z}\frac{(1+z)^{wtu+m}}{z^{m+i+1}}Y(u,z)v\\
& &+\sum\limits_{i=0}^{n}(m+i+1)(-1)^{i}{ m+i\choose i}{\rm
Res}_{z}\frac{(1+z)^{wtu+m}}{z^{m+i+2}}Y(u,z)v
\end{eqnarray*}

Thus
\begin{eqnarray*}
& &(L(-1)u+L(0)u)\bar\ast_{m}^{n}v\\
& &=\sum\limits_{i=0}^{n}i(-1)^{i}{ m+i\choose i}{\rm
Res}_{z}\frac{(1+z)^{wtu+m}}{z^{m+i+1}}Y(u,z)v\\
& &+\sum\limits_{i=0}^{n}(m+i+1)(-1)^{i}{ m+i\choose i}{\rm
Res}_{z}\frac{(1+z)^{wtu+m}}{z^{m+i+2}}Y(u,z)v\\
& &=(n+m+1)(-1)^n{ n+m\choose m}{\rm
Res}_{z}\frac{(1+z)^{wtu+m}}{z^{n+m+2}}Y(u,z)v\\
& &=(n+m+1)(-1)^n{n+m\choose m} u\circ_m^nv.
\end{eqnarray*}
Note that $L(-1)u+L(0)u\in O_n(V).$ The proof is complete. \qed

The next lemma is motivated by the commutator relation of vertex
operators and will relate the two products $u\bar\ast_{m}^{n}v$
and $u \ast_{m}^{n}v.$
\begin{lem}\label{l2.3} If $u,v\in V,$  $p_1,p_2,m\in \Z_+$  with $p_1+p_2-m\geq 0,$ then
$$u{\ast}_{m,p_{1}}^{p_{1}+p_{2}-m}v-v\ast_{m,p_{2}}^{p_{1}+p_{2}-m}u-{\rm Res}_{z}(1+z)^{wtu-1+m-p_{2}}Y(u,z)v\in
O'_{p_{1}+p_{2}-m,m}(V).$$
\end{lem}

\pf From the  definition of $O'_{p_{1}+p_{2}-m,m}(V)$, one can
easily verify that
$$Y(v,z)u\equiv(1+z)^{-wtu-wtv-2m+p_{1}+p_{2}}Y(u,\frac{-z}{1+z})v$$
modulo $O'_{p_{1}+p_{2}-m,m}(V)$ (cf. \cite{Z} and \cite{DLM2}). Hence
\begin{eqnarray*}
&
&v\ast_{m,p_{2}}^{p_{1}+p_{2}-m}u=\sum\limits_{i=0}^{p_{2}}(-1)^{i}{
p_{1}+i\choose i}{\rm
Res}_{z}\frac{(1+z)^{wtv+m}}{z^{p_{1}+i+1}}Y(v,z)u\\
& &\ \ \ \equiv\sum\limits_{i=0}^{p_{2}}(-1)^{i}{ p_{1}+i\choose
i}{\rm
Res}_{z}\frac{(1+z)^{wtv+m}}{z^{p_{1}+i+1}}(1+z)^{-wtu-wtv-2m+p_{1}+p_{2}}Y(u,\frac{-z}{1+z})v
\\
& & \ \ \ \ \ ({{\rm mod}}\ O'_{p_{1}+p_{2}-m,m}(V))\\
 & &\ \ \
=\sum\limits_{i=0}^{p_{2}}(-1)^{p_{1}}{ p_{1}+i\choose i}{\rm
Res}_{z}\frac{(1+z)^{wtu+i-1+m-p_{2}}}{z^{p_{1}+i+1}}Y(u,z)v.
\end{eqnarray*}
Recall the definition of $u{\ast}_{m,p_{1}}^{p_{1}+p_{2}-m}v.$
Then
$$u{\ast}_{m,p_{1}}^{p_{1}+p_{2}-m}v-v\ast_{m,p_{2}}^{p_{1}+p_{2}-m}u\equiv {\rm
Res}_{z}A_{p_{1},p_{2}}(z)(1+z)^{wtu-1+m-p_{2}}Y(u,z)v,$$ where
$$
A_{p_{1},p_{2}}(z)=\sum\limits_{i=0}^{p_{1}}(-1)^{i}{
p_{2}+i\choose
i}\frac{(1+z)^{p_{2}+1}}{z^{p_{2}+i+1}}-\sum\limits_{i=0}^{p_{2}}(-1)^{p_{1}}{
p_{1}+i\choose i}\frac{(1+z)^i}{z^{p_{1}+i+1}}.$$
 The lemma now follows from Proposition 5.1 in the Appendix. \qed

The proof of the following lemma is fairly standard (cf.
\cite{DLM3} and \cite{Z}).
\begin{lem}\label{l2.4} For homogeneous  $u,v\in V$, and integers
 $k\geq s\geq 0$,
 $$
{\rm Res}_{z}\frac{(1+z)^{wtu+m+s}}{z^{n+m+2+k}}Y(u,z)v\in
O'_{n,m}(V).$$
\end{lem}

\begin{lem}\label{l2.5} We have
$V\bar{\ast}_{m}^{n}O'_{n,m}(V)\subseteq O'_{n,m}(V),$
$O'_{n,m}(V)\ast_{m}^{n}V\subseteq O'_{n,m}(V)$.
\end{lem}

\pf For homogeneous $u,v\in V$ and $w\in V$,
\begin{eqnarray*}
& &\ \ \ \ \ u\bar{\ast}_{m}^{n}(v\circ_{m}^{n}w)\\
& & \equiv\sum\limits_{i=0}^{n}(-1)^{i}{ m+i\choose i}{\rm
Res}_{z_{1}}\frac{(1+z_{1})^{wtv+m}}{z_{1}^{m+i+1}}Y(u,z_{1}){\rm
Res}_{z_{2}}\frac{(1+z_{2})^{wtv+m}}{z_{2}^{n+m+2}}Y(v,z_{2})w\\
& &\ \ \ \ -\sum\limits_{i=0}^{n}(-1)^{i}{ m+i\choose i}{\rm
Res}_{z_{2}}\frac{(1+z_{2})^{wtv+m}}{z_{2}^{n+m+2}}Y(v,z_{2}){\rm
Res}_{z_{1}}\frac{(1+z_{1})^{wtv+m}}{z_{1}^{m+i+1}}Y(u,z_{1})w\\
& &=\sum\limits_{i=0}^{n}(-1)^{i}{ m+i\choose
i}\sum\limits_{j\geq0}{ wtu+m\choose j}\sum\limits_{k=0}^{\infty}
{-m-i-1\choose k}\\
& &\ \ \ \ \cdot {\rm Res}_{z_{2}}{\rm
Res}_{z_{1}-z_{2}}\frac{(1+z_{2})^{wtu+wtv+2m-j}(z_{1}-z_{2})^{j+k}}{z_{2}^{2m+n+i+3+k}}
Y(Y(u,z_{1}-z_{2})v,z_{2})w\\
& &=\sum\limits_{i=0}^{n}(-1)^{i}{ m+i\choose
i}\sum\limits_{j\geq0}{ wtu+m\choose j}\sum\limits_{k=0}^{\infty}
{-m-i-1\choose k}\\
& &\ \ \ \  \ \cdot{\rm
Res}_{z_{2}}\frac{(1+z_{2})^{wtu+wtv+2m-j}}{z_{2}^{2m+n+i+3+k}}Y(u_{j+k}v,z_{2})w.
\end{eqnarray*}
Note that the weight of $u_{j+k}v$ is $wtu+wtv-j-k-1.$ By Lemma
\ref{l2.4} we see that $u\bar{\ast}_{m}^{n}(v\circ_{m}^{n}w)$ lies
in $O'_{n,m}(V).$

By Lemma \ref{l2.3}, we have
\begin{eqnarray*}
& &\ \ \ \ u\bar{\ast}_{m}^{n}(v\circ_{m}^{n}w)-(v\circ_{m}^{n}w){\ast}_{m}^{n}u\\
& &\equiv{\rm Res}_{z}(1+z)^{wtu-1}Y(u,z)(v\circ_{m}^{n}w)\\
& &={\rm Res}_{z_{1}}(1+z_{1})^{wtu-1}{\rm
Res}_{z_{2}}\frac{(1+z_{2})^{wtv+m}}{z_{2}^{n+m+2}}Y(u,z_{1})Y(v,z_{2})w\\
& &\equiv\sum\limits_{j\geq 0}{ wtu-1\choose j}{\rm
Res}_{z_{2}}{\rm
Res}_{z_{1}-z_{2}}\frac{(1+z_{2})^{wtu+wtv+m-1-j}(z_{1}-z_{2})^{j}}{z_{2}^{n+m+2}}\\
& &\ \ \ \ \cdot Y(Y(u,z_{1}-z_{2})v,z_{2})w\\
& &=\sum\limits_{j\geq 0}{ wtu-1\choose j}{\rm
Res}_{z_{2}}\frac{(1+z_{2})^{wtu+wtv+m-1-j}}{z_{2}^{n+m+2}}Y(u_jv,z_{2})w\in
O'_{n,m}(V)
\end{eqnarray*}
which is a vector in $O'_{n,m}(V)$ by the definition of
$O'_{n,m}(V)$. As a result, $(v\circ_{m}^{n}w){\ast}_{m}^{n}u\in
O'_{n,m}(V).$

Next we deal
with $L(-1)u+(L(0)+m-n)u\in O'_{n,m}(V).$ As before we assume that $u$
is homogeneous. Then
\begin{eqnarray*}
& &\ \ \ \ (L(-1)u+(L(0)+m-n)u){\ast}_{m}^{n}v\\
& &=\sum\limits_{i=0}^{m}(-1)^{i}{ n+i\choose i}{\rm
Res}_{z}\frac{(1+z)^{wtu+m+1}}{z^{n+i+1}}Y(L(-1)u,z)v\\
& &\ \ \ \ +\sum\limits_{i=0}^{m}(-1)^{i} {n+i\choose i}
(wtu+m-n){\rm
Res}_{z}\frac{(1+z)^{wtu+m}}{z^{n+i+1}}Y(u,z)v\\
& &=\sum\limits_{i=0}^{m}(-1)^{i+1}{ n+i\choose i} (wtu+m+1){\rm
Res}_{z}\frac{(1+z)^{wtu+m}}{z^{n+i+1}}Y(u,z)v\\
& &\ \ \ \ +\sum\limits_{i=0}^{m}(-1)^{i}{ n+i\choose i}
(n+i+1){\rm
Res}_{z}\frac{(1+z)^{wtu+m+1}}{z^{n+i+2}}Y(u,z)v\\
& &\ \ \  \ +\sum\limits_{i=0}^{m}(-1)^{i}{ n+i\choose i}
(wtu+m-n){\rm Res}_{z}\frac{(1+z)^{wtu+m}}{z^{n+i+1}}Y(u,z)v\\
& &=\sum\limits_{i=0}^{m}(-1)^{i+1}{ n+i\choose i} (n+1){\rm
Res}_{z}\frac{(1+z)^{wtu+m}}{z^{n+i+1}}Y(u,z)v\\
& &\ \ \ \ +\sum\limits_{i=0}^{m}(-1)^{i}{ n+i\choose i}
(n+i+1){\rm
Res}_{z}\frac{(1+z)^{wtu+m}}{z^{n+i+1}}Y(u,z)v\\
& &\ \ \  \ +\sum\limits_{i=0}^{m}(-1)^{i}{ n+i\choose i}
(n+i+1){\rm
Res}_{z}\frac{(1+z)^{wtu+m}}{z^{n+i+2}}Y(u,z)v.
\end{eqnarray*}
It is easy to show that the last expression is equal to $(-1)^m(n+m+1)
{m+n\choose m}u\circ_n v.$ By Lemma \ref{l2.2}, $(L(-1)u+(L(0)+m-n)u){\ast}_{m}^{n}v$ belongs to $O'_{n,m}(V).$

Finally for $v\bar{\ast}_{m}^{n}(L(-1)u+(L(0)+m-n)u)$
we use  Lemma \ref{l2.3} to yield
\begin{eqnarray*}
& &\ \ \ \ v\bar{\ast}_{m}^{n}(L(-1)u+(L(0)+m-n)u)-(L(-1)u+(L(0)+m-n)u){\ast}_{m}^{n}v\\
& &={\rm Res}_{z}(1+z)^{wtv-1}Y(v,z)(L(-1)u+(L(0)+m-n)u)\\
& &=\sum\limits_{i\geq0}{ wtv-1\choose
i}v_iL(-1)u+(wtu+m-n)\sum\limits_{i\geq 0}{ wtv-1\choose i}v_iu\\
& &=L(-1)\sum\limits_{i\geq 0}{ wtv-1\choose
i}v_iu+\sum\limits_{i\geq 0}{ wtv-1\choose i}iv_{i-1}u\\
& &\ \ \ \ +(wtu+m-n)\sum\limits_{i\geq 0}
{ wtv-1\choose i}v_iu\\
& &=L(-1)\sum\limits_{i\geq 0}{ wtv-1\choose i}
v_iu+\sum\limits_{i\geq 0}{ wtv-1\choose i+1}(i+1)v_iu\\
& &\ \ \ \ +(wtu+m-n)\sum\limits_{i\geq 0}{ wtv-1\choose i}
v_iu\\
& &=\sum\limits_{i\geq 0} { wtv-1\choose
i}(L(-1)+wtv-i-1+wtu+m-n)v_iu\\
& &= \sum\limits_{i\geq 0} { wtv-1\choose i}(L(-1)v_iu+ L(0)v_iu+
(m-n)v_iu)
\end{eqnarray*}
which is in $O'_{n,m}(V).$ So $v\bar{\ast}_{m}^{n}(L(-1)u+(L(0)+m-n)u)\in
O'_{n,m}(V),$ as desired. \qed

We should remind the reader that our goal is to construct an
$A_n(V)$-$A_m(V)$-bimodule $A_{n,m}(V)$ with the left action $\bar
*_{m}^n$ of $A_n(V)$ and the right action $*_m^n$ of $A_m(V).$ The
following lemma claims that the left action $\bar *_{m}^n$ and the
right action $*_m^n$ commute. On the other hand, we do not need to
prove this lemma as a bigger subspace
$O_{n,m}(V)$ of $V$ containing $O'_{n,m}(V)$
will be modulo out. In fact,
$(a\bar{\ast}_{m}^{n}b)\ast_{m}^{n}c-a\bar{\ast}_{m}^{n}(b\ast_{m}^{n}c)$
is an element of $O_{n,m}(V)$ (see Lemma \ref{l2.6} below). But eventually we expect to prove
that $O_{m,n}(V)$ and $O'_{n,m}(V)$ are the same although we
cannot achieve this in the paper.

\begin{lem}\label{l2.6} We have $
(a\bar{\ast}_{m}^{n}b)\ast_{m}^{n}c-a\bar{\ast}_{m}^{n}(b\ast_{m}^{n}c)$
lies in $O'_{n,m}(V)$ for homogeneous $a,b,c\in V.$
\end{lem}

\pf The proof of this lemma is similar to that of Theorem 2.4 of \cite{DLM3}.
In fact, if $m=n,$ the lemma is exactly the associativity of product
$*_n$ in $A_n(V).$

A straightforward calculation using Lemma \ref{l2.4} gives:
\begin{eqnarray*}
& & \ \ \ \ (a\bar{\ast}_{m}^{n}b)\ast_{m}^{n}c\\
& &=\sum\limits_{k=0}^{m}(-1)^{k} {n+k\choose
k}\sum\limits_{i=0}^{n}(-1)^{i} {m+i\choose i}\sum\limits_{j\geq
0}
{wta+m\choose j}\\
& &\ \ \ \ \cdot{\rm
Res}_{z}\frac{(1+z)^{wta+wtb+2m-j+i}}{z^{n+k+1}}Y(a_{j-m-i-1}b,z)c
\\
& &=\sum\limits_{k=0}^{m}(-1)^{k}{n+k\choose k}
\sum\limits_{i=0}^{n}(-1)^{i}{m+i\choose i}
\sum\limits_{j\geq 0}{wta+m\choose j}\\
& &\ \ \ \ \cdot{\rm Res}_{z_{2}}{\rm
Res}_{z_{1}-z_{2}}\frac{(1+z_{2})^{wta+wtb+2m-j+i}(z_{1}-z_{2})^{j-m-i-1}}{z_{2}^{n+k+1}}Y(Y(a,z_{1}-z_{2})b,z_{2})c\\
& &=\sum\limits_{k=0}^{m}(-1)^{k}{n+k\choose
k}\sum\limits_{i=0}^{n}(-1)^{i}{m+i\choose i}{\rm Res}_{z_{2}}{\rm
Res}_{z_{1}-z_{2}}\\
& &\ \ \ \ \cdot\frac{(1+z_{1})^{wta+m}(1+z_{2})^{wtb+m+i}(z_{1}-z_{2})^{-m-i-1}}{z_{2}^{n+k+1}}
Y(Y(a,z_{1}-z_{2})b,z_{2})c\\
& &=\sum\limits_{k=0}^{m}(-1)^{k}{n+k\choose
k}\sum\limits_{i=0}^{n}(-1)^{i} {m+i\choose i}\sum\limits_{j\geq
0}
{-m-i-1\choose j}\\
& &\ \ \ \ \cdot{\rm Res}_{z_{1}}{\rm Res}_{z_{2}}
\frac{(1+z_{1})^{wta+m}(1+z_{2})^{wtb+m+i}(-z_{2})^{j}}{z_{1}^{m+i+1+j}z_{2}^{n+k+1}}Y(a,z_{1})Y(b,z_{2})c\\
& &\ \ \ \ -\sum\limits_{k=0}^{m}(-1)^{k} {n+k\choose
k}\sum\limits_{i=0}^{n}(-1)^{i} {m+i\choose i}\sum\limits_{j\geq
0}
{-m-i-1\choose j}\\
& &\ \ \ \ \cdot{\rm Res}_{z_{2}}{\rm Res}_{z_{1}}
\frac{(1+z_{1})^{wta+m}(1+z_{2})^{wtb+m+i}z_{1}^{j}}{(-z_{2})^{m+i+1+j}z_{2}^{n+k+1}}Y(b,z_{2})Y(a,z_{1})c\\\
& & \equiv
a\bar{\ast}_{m}^{n}(b\ast_{m}^{n}c)+\sum\limits_{k=0}^{m}(-1)^{k}
{n+k\choose k}\sum\limits_{i=0}^{n}(-1)^{i}
{m+i\choose i}\\
& &\ \ \ \ \cdot{\rm Res}_{z_{1}}{\rm
Res}_{z_{2}}\left[\sum\limits_{j=0}^{n-i}(-1)^{j}
{-m-i-1\choose j}\sum\limits_{l=0}^{i}{ i\choose l}\frac{z_{2}^{j+l}}{z_{1}^{j+i}}-\frac{1}{z_{1}^{i}}\right]\\
& &\ \ \ \
\cdot\frac{(1+z_{1})^{wta+m}(1+z_{2})^{wtb+m}}{z_{1}^{m+1}z_{2}^{n+k+1}}Y(a,z_{1})Y(b,z_{2})c.
\end{eqnarray*}
The lemma then follows from Proposition \ref{p5.2} in the Appendix. \qed

In order to construct $A_{n,m}(V)$ we need to introduce more
subspaces of $V.$ Let $O''_{n,m}(V)$ be the linear span of
$u\ast_{m,p_{3}}^{n}((a\ast_{p_{1},p_{2}}^{p_{3}}b){\ast}_{m,p_{1}}^{p_{3}}c-
a{\ast}_{m,p_{2}}^{p_{3}}(b{\ast}_{m,p_{1}}^{p_{2}}c)),$
 for $a,b,c,u\in V, p_{1},p_{2},p_{3}\in{\mathbb
Z}_{+}$, and $O'''_{n,m}(V)=\sum_{p\in{\mathbb
Z}_{+}}(V\ast_{p}^{n}O_{p}(V))\ast_{m,p}^{n}V.$ Set
$$
O_{n,m}(V)=O'_{n,m}(V)+O''_{n,m}(V)+O'''_{n,m}(V).
$$

\begin{lem} For $p_{1},p_{2},m,n\in{\mathbb Z}_{+}$, we have
$(V\ast_{p_{1},p_{2}}^{n}O'_{p_{2},p_{1}}(V))\ast_{m,p_{1}}^{n}V\subseteq
O_{n,m}(V)$.
\end{lem}

\pf  We first prove that
$(L(-1)u+(L(0)+p_{1}-p_{2})u)\ast_{m,p_{1}}^{p_{2}}v\in
O_{p_{2}}(V)\ast_{m,p_{2}}^{p_{2}}V$, for homogeneous $u\in V$ and
$v\in V$. In fact,
\begin{eqnarray*}
& &\ \ \ \ (L(-1)u+(L(0)+p_{1}-p_{2})u){\ast}_{m,p_{1}}^{p_{2}}v\\
& &=\sum\limits_{i=0}^{p_{1}}(-1)^{i}{ m+p_{2}-p_{1}+i\choose
i}{\rm
Res}_{z}\frac{(1+z)^{wtu+m+1}}{z^{m+p_{2}-p_{1}+i+1}}Y(L(-1)u,z)v\\
& &\ \ \ \ +\sum\limits_{i=0}^{p_{1}}(-1)^{i}
{m+p_{2}-p_{1}+i\choose i} (wtu+p_{1}-p_{2}){\rm
Res}_{z}\frac{(1+z)^{wtu+m}}{z^{m+p_{2}-p_{1}+i+1}}Y(u,z)v\\
& &=\sum\limits_{i=0}^{p_{1}}(-1)^{i+1}{ m+p_{2}-p_{1}+i\choose i}
(wtu+m+1){\rm
Res}_{z}\frac{(1+z)^{wtu+m}}{z^{m+p_{2}-p_{1}+i+1}}Y(u,z)v\\
& &\ \ \ \ +\sum\limits_{i=0}^{p_{1}}(-1)^{i}{
m+p_{2}-p_{1}+i\choose i} (m+p_{2}-p_{1}+i+1){\rm
Res}_{z}\frac{(1+z)^{wtu+m+1}}{z^{m+p_{2}-p_{1}+i+2}}Y(u,z)v\\
& &\ \ \  \ +\sum\limits_{i=0}^{p_{1}}(-1)^{i}{
m+p_{2}-p_{1}+i\choose i} (wtu+p_{1}-p_{2}){\rm
Res}_{z}\frac{(1+z)^{wtu+m}}{z^{m+p_{2}-p_{1}+i+1}}Y(u,z)v\\
& &=(-1)^{p_{1}}{m+p_{2}\choose p_{1}}(m+p_{2}+1){\rm
Res}_{z}\frac{(1+z)^{wtu+m}}{z^{m+p_{2}+2}}Y(u,z)v.
\end{eqnarray*}
So by Lemma 2.2,
$(L(-1)u+(L(0)+p_{1}-p_{2})u)\ast_{m,p_{1}}^{p_{2}}v\in
O_{p_{2}}(V)\ast_{m,p_{2}}^{p_{2}}V$. Hence by the definitions of
$O'_{n,m}(V)$ and $O_{n,m}(V),$ and Lemma 2.2, we have
\begin{eqnarray*}
& &(V\ast_{p_{1},p_{2}}^{n}O'_{p_{2},p_{1}}(V))\ast_{m,p_{1}}^{n}V\subseteq
V\ast_{m,p_{2}}^{n}(O'_{p_{2},p_{1}}(V)\ast_{m,p_{1}}^{p_{2}}V)+O_{n,m}(V)
\\
& &\subseteq
V\ast_{m,p_{2}}^{n}((O_{p_{2}}(V)\ast_{p_{1},p_{2}}^{p_{2}}V)\ast_{m,p_{1}}^{p_{2}}V+
O_{p_{2}}(V)\ast_{m,p_{2}}^{p_{2}}V)+O_{n,m}(V)\\
& &\subseteq
V\ast_{m,p_{2}}^{n}(O_{p_{2}}(V)\ast_{m,p_{2}}^{p_{2}}(V\ast_{m,p_{1}}^{p_{2}}V))+
(V\ast_{p_{2}}^{n}O_{p_{2}}(V))\ast_{m,p_{2}}^{n}V+O_{n,m}(V)\\
& &
\subseteq
(V\ast_{p_{2}}^{n}O_{p_{2}}(V))\ast_{m,p_{2}}^{n}(V\ast_{m,p_{1}}^{p_{2}}V)+O_{n,m}(V)\\
& &\subseteq
O_{n,m}(V),
\end{eqnarray*}
as required. \qed

\begin{lem}\label{l2.8} For
any $m,n,p\in{\mathbb Z}_{+},$ we have  $V{\ast}_{m,p}^{n}O_{p,m}(V)\subseteq
O_{n,m}(V)$, $O_{n,p}(V)\ast_{m,p}^{n}V\subseteq O_{n,m}(V).$ In particular,
$V\bar{\ast}_{m}^{n}O_{n,m}(V)\subseteq O_{n,m}(V)$,
$O_{n,m}(V)\ast_{m}^{n}V\subseteq O_{n,m}(V)$.
\end{lem}

\pf  By   the definition of $O_{n,m}(V)$ and Lemma 2.7, it
suffices to prove that
\begin{equation}\label{2.7}
V{\ast}_{m,p}^{n}((V\ast_{p_{1}}^{p}O_{p_{1}}(V)){\ast}_{m,p_{1}}^{p}V+O''_{p,m}(V))\subseteq
O_{n,m}(V)
\end{equation}
and
\begin{equation}\label{2.8}
((V\ast_{p_{1}}^{n}O_{p_{1}}(V)){\ast}_{p,p_{1}}^{n}V+O''_{n,p}(V))\ast_{m,p}^{n}V\subseteq
O_{n,m}(V)
\end{equation}
for $p_1,p\in\Z_+.$

 We first prove (2.1). It is clear that
\begin{eqnarray*}
& &V{\ast}_{m,p}^{n}((V\ast_{p_{1}}^{p}O_{p_{1}}(V)){\ast}_{m,p_{1}}^{p}V)\subseteq
V{\ast}_{m,p}^{n}(V\ast_{m,p_{1}}^{p}(O_{p_{1}}(V){\ast}_{m,p_{1}}^{p_{1}}V))+O_{n,m}(V)
 \\
& &\subseteq
(V{\ast}_{p_{1},p}^{n}V)\ast_{m,p_{1}}^{n}(O_{p_{1}}(V){\ast}_{m,p_{1}}^{p_{1}}V)+O_{n,m}(V)\\
& &\subseteq
((V{\ast}_{p_{1},p}^{n}V)\ast_{p_{1}}^{n}O_{p_{1}}(V)){\ast}_{m,p_{1}}^{n}V+O_{n,m}(V)\subseteq
O_{n,m}(V).
\end{eqnarray*}
 It remains to prove that
$V{\ast}_{m,p}^{n}O''_{p,m}(V)\subseteq O_{n,m}(V).$ With $u={\bf
1}$ in the definition of $O''_{n,m}(V)$ we have
$(a\ast_{p_{1},p_{2}}^{n}b){\ast}_{m,p_{1}}^{n}c-
a{\ast}_{m,p_{2}}^{n}(b{\ast}_{m,p_{1}}^{p_{2}}c)\in O_{n,m}(V).$
Thus
\begin{eqnarray*}
& & \ \ \ \ v{\ast}_{m,p}^{n}
(u\ast_{m,p_{3}}^{p}((a\ast_{p_{1},p_{2}}^{p_{3}}b){\ast}_{m,p_{1}}^{p_{3}}c-
a{\ast}_{m,p_{2}}^{p_{3}}(b{\ast}_{m,p_{1}}^{p_{2}}c)))\\
& &\equiv
(v{\ast}_{p_{3},p}^{n}u)\ast_{m,p_{3}}^{n}((a\ast_{p_{1},p_{2}}^{p_{3}}b){\ast}_{m,p_{1}}^{p_{3}}c-
a{\ast}_{m,p_{2}}^{p_{3}}(b{\ast}_{m,p_{1}}^{p_{2}}c))\\& & \equiv
0 \ ({\rm mod}\  O_{n,m}(V)).
\end{eqnarray*}
So (2.1) is true.

For (2.2), it is easy to see  that
$((V\ast_{p_{1}}^{n}O_{p_{1}}(V)){\ast}_{p,p_{1}}^{n}V)\ast_{m,p}^{n}V\subseteq
(V\ast_{p_{1}}^{n}O_{p_{1}}(V)){\ast}_{m,p_{1}}^{n}(V\ast_{m,p}^{p_{1}}V)+O_{n,m}(V)\subseteq
O_{n,m}(V)$. Let $a,b,c,u,v\in V$ and
$p_{1},p_{2},p_{3},p\in{\mathbb Z}_{+}$, then
\begin{eqnarray*} &
&
(u\ast_{p,p_{3}}^{n}((a\ast_{p_{1},p_{2}}^{p_{3}}b){\ast}_{p,p_{1}}^{p_{3}}c-
a{\ast}_{p,p_{2}}^{p_{3}}(b{\ast}_{p,p_{1}}^{p_{2}}c)))\ast_{m,p}^{n}v\\
& & \equiv
u\ast_{m,p_{3}}^{n}(((a\ast_{p_{1},p_{2}}^{p_{3}}b){\ast}_{p,p_{1}}^{p_{3}}c-
a{\ast}_{p,p_{2}}^{p_{3}}(b{\ast}_{p,p_{1}}^{p_{2}}c))\ast_{m,p}^{p_{3}}v)\\
& &\equiv
u\ast_{m,p_{3}}^{n}((a\ast_{p_{1},p_{2}}^{p_{3}}b)\ast_{m,p_{1}}^{p_{3}}(c\ast_{m,p}^{p_{1}}v))
-u\ast_{m,p_{3}}^{n}(a\ast_{m,p_{2}}^{p_{3}}((b\ast_{p,p_{1}}^{p_{2}}c)\ast_{m,p}^{p_{2}}v))\\
& & \equiv
u\ast_{m,p_{3}}^{n}(a\ast_{m,p_{2}}^{p_{3}}(b\ast_{m,p_{1}}^{p_{2}}(c\ast_{m,p}^{p_{1}}v)))-
(u\ast_{p_{2},p_{3}}^{n}a)\ast_{m,p_{2}}^{n}((b\ast_{p,p_{1}}^{p_{2}}c)\ast_{m,p}^{p_{2}}v)\\
& & \equiv
(u\ast_{p_{2},p_{3}}^{n}a)\ast_{m,p_{2}}^{n}(b\ast_{m,p_{1}}^{p_{2}}(c\ast_{m,p}^{p_{1}}v))
-(u\ast_{p_{2},p_{3}}^{n}a)\ast_{m,p_{2}}^{n}((b\ast_{p,p_{1}}^{p_{2}}c)\ast_{m,p}^{p_{2}}v)\\
 & & \equiv 0 \ ({\rm mod} \ O_{n,m}(V)).
\end{eqnarray*}
The lemma is proved. \qed
\\

 We now define
$$A_{n,m}(V)=V/O_{n,m}(V).$$
The reason for this definition will become clear from the
representation theory of $V$ discussed later. The following is the
first main theorem in this paper.

\begin{theorem}\label{t2.8} Let $V$ be a vertex operator algebra and $m,n$ nonnegative integers.
Then  $A_{n,m}(V)$ is an $A_n(V)$-$A_m(V)$-bimodule such that the
left and right actions of $A_n(V)$ and $A_m(V)$ are given by
$\bar*_m^n$ and $*^n_m.$
\end{theorem}

\pf First, both actions are well defined from the definition of
$O_{n,m}(V)$ and Lemma \ref{l2.8}. The left $A_n(V)$-module and
right $A_m(V)$-module structures then follow from the fact that
$O_{n,m}''(V)$ is a subspace of $O_{n,m}(V).$ The commutativity of
two actions proved in Lemma \ref{l2.6} asserts that   $A_{n,m}(V)$
is an $A_n(V)$-$A_m(V)$-bimodule. \qed

\begin{remark} We will prove in Section 4 that
if $m=n,$ the $A_{n,n}(V)$ defined here is the same as $A_n(V)$
discussed before. In particular, $O_{n,n}(V)$ and $O'_{n,n}(V)$
coincide. In other words, $O''_{n,n}(V),$ $O'''_{n,n}(V)$ are
subspaces of $O'_{n,n}(V).$ We suspect that this is true in
general. That is, $O_{n,m}(V)$ and its subspace $O'_{n,m}(V)$ are
equal. It seems that this is a very difficult problem and we
cannot find a proof for this in this paper.
\end{remark}

\section{Properties of $A_{n.m}(V)$}
\def\theequation{3.\arabic{equation}}
\setcounter{equation}{0} In this section we will discuss some
important properties of $A_{n.m}(V)$ such as isomorphism between
$A_{n,m}(V)$ and $A_{m,n}(V),$ relations between $A_{n,m}(V)$ and
$A_{l,k}(V)$ and tensor products. Some of  these properties will
be interpreted in terms representation theory in the next section.

First we establish the isomorphism between $A_{n,m}(V)$ and
$A_{m,n}(V)$ as $A_n(V)$-$A_m(V)$-bimodules. To achieve this we
need to define new actions of $A_n(V)$ and $A_m(V)$ on
$A_{m,n}(V)$ so that $A_{m,n}(V)$ becomes an
$A_n(V)$-$A_m(V)$-bimodule. Recall from  \cite{Z} the linear map
$\phi:V\to V$ such that  $\phi(v)=e^{L(1)}(-1)^{L(0)}v$ for $v\in V.$
Then $\phi$ induces an anti involution on $A_n(V)$ \cite{DLM3}
(also see \cite{Z},\cite{DLM2}). We
also use $\phi$ at the present situation
to define an isomorphism between $A_{n,m}(V)$ and
$A_{m,n}(V).$
\begin{lem}\label{l3.1} For $u,v\in V$ define
$$u\bar\cdot_m^nv=v*_n^m\phi(u), \ \ u\cdot_m^n v=\phi(v)\bar *^m_nu.$$
Then $A_{m,n}(V)$ becomes an $A_n(V)$-$A_m(V)$-bimodule under the
left action $\bar\cdot_m^n$ by $A_n(V)$ and the right action
$\cdot_m^n$ by $A_m(V).$
\end{lem}

\pf Since $\phi(O_m(V))\subset O_m(V)$ for any $m$ (see
\cite{DLM3}), we immediately see that both actions are well
defined. The rest follows from Theorem \ref{t2.8} and the fact
that $\phi$ is an anti involution of $A_m(V)$ for any $m.$ \qed

\begin{prop} The linear map $ \phi$ induces an isomorphism of
$A_{n}(V)$-$A_m(V)$-bimodules from $A_{n,m}(V)$ to $A_{m,n}(V)$,
where the actions of $A_{n}(V)$ and $A_m(V)$ on $A_{n,m}(V)$ are
defined in Theorem \ref{t2.8}, and the actions on $A_{m,n}(V)$ are
defined in Lemma \ref{l3.1}.
\end{prop}

\pf  We first prove that
 $\phi(u{\ast}_{m,p}^{n}v)\equiv \phi(v){\ast}_{n,p}^{m}\phi(u)$
modulo $O'_{m,n}(V)$ for $u,v\in V$ and $p\in\Z_+.$  Recall the identities
$$(-1)^{L(0)}Y(u,z)(-1)^{L(0)}=Y((-1)^{L(0)}u,-z)$$
$$e^{L(1)}Y(u,z)e^{-L(1)}=Y(e^{(1-z)L(1)}(1-z)^{-2L(0)}u,\frac{z}{1-z})$$
from \cite{FLM}. We have the following computation with the help
from Proposition \ref{p5.1} in Appendix:
\begin{eqnarray*}
&&\phi(u{\ast}_{m,p}^{n}v)=\phi\left(\sum\limits_{i=0}^{p}(-1)^{i}
{m+n-p+i\choose i}{\rm Res}_{z}\frac{(1+z)^{wtu+m}}{z^{m+n-p+i+1}}Y(u,z)v\right)\\
& &=\sum\limits_{i=0}^{p}(-1)^{i}{m+n-p+i\choose i}{\rm
Res}_{z}\frac{(1+z)^{wtu+m}}{z^{m+n-p+i+1}}e^{L(1)}Y((-1)^{L(0)}u,-z)(-1)^{L(0)}v\\
& &=\sum\limits_{i=0}^{p}(-1)^{i}{m+n-p+i\choose i}{\rm
Res}_{z}\frac{(1+z)^{wtu+m}}{
z^{m+n-p+i+1}}\\
& &\ \ \ \ \cdot
Y(e^{(1+z)L(1)}(1+z)^{-2L(0)}(-1)^{L(0)}u,\frac{-z}{1+z})e^{L(1)}(-1)^{L(0)}v
\\
& &=\sum\limits_{i=0}^{p}(-1)^{wtu+m+n-p}{m+n-p+i\choose i}{\rm
Res}_{z}\frac{(1+z)^{wtu+n-p+i-1}}{
z^{m+n-p+i+1}}\\
& & \ \ \ \ \ \ \cdot Y(e^{(1+z)^{-1}L(1)}u,z)e^{L(1)}(-1)^{L(0)}v\\
& &=\sum\limits_{j=0}^{\infty}\frac{1}{
j!}\sum\limits_{i=0}^{p}(-1)^{wtu+m+n-p}{m+n-p+i\choose i}{\rm
Res}_{z}\frac{(1+z)^{wtu+n-p-j+i-1}}{
z^{m+n-p+i+1}}\\
& & \ \ \ \ \ \ \cdot Y(L(1)^{j}u,z)e^{L(1)}(-1)^{L(0)}v\\
& &=\sum\limits_{j=0}^{\infty}\frac{(-1)^{wtu}}{
j!}\sum\limits_{i=0}^{m+n-p}(-1)^{i}\left(\begin{array}{c}
p+i\\i\end{array} \right){\rm Res}_{z}\frac{(1+z)^{wtu-j+n}}{
z^{p+i+1}}Y(L(1)^{j}u,z)e^{L(1)}(-1)^{L(0)}v\\
& &\ \ \ \ -\sum\limits_{j=0}^{\infty}\frac{(-1)^{wtu}}{j!}{\rm
Res}_{z}(1+z)^{wtu-j-1+n-p}Y(L(1)^{j}u,z)e^{L(1)}(-1)^{L(0)}v\\
& &\equiv \phi(v){\ast}_{n,p}^{m}\phi(u)\ ({\rm mod}\
O'_{m,n}(V)),
\end{eqnarray*}
where we have used Proposition \ref{p5.1} and Lemma \ref{l2.3} in
the last two steps. In particular,
$\phi(u\bar{\ast}_{m}^{n}v)\equiv \phi(v)\ast_{n}^{m}\phi(u)$
modulo $O'_{n,m}(V)$ and $\phi(u{\ast}_{m}^{n}v)\equiv
\phi(v)\bar{\ast}_{n}^{m}\phi(u)$ modulo $O'_{m,n}(V)$ for $u,v\in
V.$

We next prove that $\phi(O_{n,m}(V))\subset O_{m,n}(V)$. Take
$v={\rm Res}_{z}\frac{(1+z)^{wta+m}}{z^{n+m+2}}Y(a,z)b\in
O'_{n,m}(V).$ Then
\begin{eqnarray*}
& &\phi(v)={\rm Res}_{z}e^{L(1)}(-1)^{L(0)}\frac{(1+z)^{wta+m}}{
z^{n+m+2}}Y(a,z)b\\
& &={\rm Res}_{z}e^{L(1)}\frac{(1+z)^{wta+m}}{
z^{n+m+2}}Y(e^{(1+z)L(1)}(1+z)^{-2L(0)}(-1)^{L(0)}a,\frac{-z}{1+z})e^{L(1)}(-1)^{L(0)}b\\
& &={\rm Res}_{z}(-1)^{wta+m+n+1}\frac{(1+z)^{wta+n}}{
z^{n+m+2}}Y(e^{\frac{1}{{1+z}}L(1)}a,z)e^{L(1)}(-1)^{L(0)}b\\
& &={\rm
Res}_{z}(-1)^{wta+m+n+1}\sum\limits_{j=0}^{\infty}
\frac{1}{j!}\frac{(1+z)^{wta+n-j}}{z^{n+m+2}}Y(L(1)^{j}a,z)e^{L(1)}(-1)^{L(0)}b\in
O'_{m,n}(V).
\end{eqnarray*}

For $u\in V$, we have
\begin{eqnarray*}
& &\ \ \ \ \phi(L(-1)u+(L(0)+m-n)u)\\
& &=e^{L(1)}(-1)^{L(0)}{\rm
Res}_{z}(Y(\omega,z)u+zY(\omega,z)u)+(m-n)e^{L(1)}(-1)^{L(0)}u\\
& &=e^{L(1)}{\rm
Res}_{z}(Y(\omega,-z)+zY(\omega,-z))(-1)^{L(0)}u+(m-n)e^{L(1)}(-1)^{L(0)}u\\
& &={\rm
Res}_{z}(1+z)Y(e^{(1+z)L(1)}(1+z)^{-2L(0)}(-1)^{L(0)}\omega,{\frac{-z}{1+z}})e^{L(1)}(-1)^{L(0)}u\\
& &\ \ \ \ +(m-n)e^{L(1)}(-1)^{L(0)}u \\
& &={\rm Res}_{z}(-(1+z)^{2}+z(1+z))Y(e^{(1+z)^{-1}L(1)}\omega,z)e^{L(1)}(-1)^{L(0)}u\\
& &\ \ \ \ +(m-n)e^{L(1)}(-1)^{L(0)}u\\
& &=-(L(-1)+L(0))e^{L(1)}(-1)^{L(0)}u-(n-m)e^{L(1)}(-1)^{L(0)}u\in
O'_{m,n}(V).
\end{eqnarray*}
So $\phi(O'_{n,m}(V))\subset O'_{m,n}(V)$.

Let $u\in O'_{n,p}(V),\ w\in O'_{p,m}(V)$, $v\in V.$  Then
$\phi(u)\in O'_{p,n}(V)$ and $\phi(w)\in O'_{m,p}(V).$ From the
proof above, we have $\phi (u{\ast}_{m,p}^{n}v)\equiv
\phi(v)\ast_{n,p}^{m}\phi(u)$ and $\phi(v{\ast}_{m,p}^{n}w)\equiv
\phi(w){\ast}_{n,p}^{m}\phi(v)$ modulo $O'_{m,n}(V).$ Thus by
Lemma 2.7 and the definition of $O_{m,n}(V)$,
\begin{eqnarray*}
& &\phi((V\ast_{p_{1}}^{n}O_{p_{1}}(V))\ast_{m,p_{1}}^{n}V)
\subseteq
V\ast_{n,p_{1}}^{m}\phi(V\ast_{p_{1}}^{n}O_{p_{1}}(V))+O'_{m,n}(V)\\
& &\subseteq
V\ast_{n,p_{1}}^{m}(O_{p_{1}}(V)\ast_{n,p_{1}}^{p_{1}}V+O'_{p_{1},n}(V))+O_{m,n}(V)\\
& &\subseteq(V\ast_{p_{1}}^{m}O_{p_{1}}(V))\ast_{n,p_{1}}^{m}V+O_{m,n}(V)\\
& &\subseteq
O_{m,n}(V).
\end{eqnarray*}
That is, $\phi(O'''_{n,m}(V))\subset O_{m,n}(V).$

Finally we deal with $O''_{n,m}(V).$  For $u,a,b,c\in V$, by the
discussion above, we have
\begin{eqnarray*}
& &\ \ \ \
\phi[u{\ast}_{m,p_{3}}^{n}((a\ast_{p_{1},p_{2}}^{p_{3}}b){\ast}_{m,p_{1}}^{p_{3}}c-
a{\ast}_{m,p_{2}}^{p_{3}}(b{\ast}_{m,p_{1}}^{p_{2}}c))]\\
&
&\equiv\phi((a\ast_{p_{1},p_{2}}^{p_{3}}b){\ast}_{m,p_{1}}^{p_{3}}c-
a{\ast}_{m,p_{2}}^{p_{3}}(b{\ast}_{m,p_{1}}^{p_{2}}c))\ast_{n,p_{3}}^{m}\phi(u)
\ ({\rm mod} \ O'_{m,n}(V))\\
 & &
\equiv[\phi(c)\ast_{p_{3},p_{1}}^{m}\phi(a\ast_{p_{1},p_{2}}^{p_{3}}b)-
\phi(b{\ast}_{m,p_{1}}^{p_{2}}c)\ast_{p_{3},p_{2}}^{m}\phi(a)+x_{1}]\ast_{n,p_{3}}^{m}\phi(u),
\end{eqnarray*}
for some $x_{1}\in O'_{m,p_{3}}(V)$. Since by Lemma 2.7,
$O'_{m,p_{3}}(V)\ast_{n,p_{3}}^{m}V\subseteq O_{m,n}(V)$, we have
\begin{eqnarray*}
& &\ \ \ \
\phi[u{\ast}_{m,p_{3}}^{n}((a\ast_{p_{1},p_{2}}^{p_{3}}b){\ast}_{m,p_{1}}^{p_{3}}c
-
a{\ast}_{m,p_{2}}^{p_{3}}(b{\ast}_{m,p_{1}}^{p_{2}}c))]\\
 & &\equiv
[\phi(c)\ast_{p_{3},p_{1}}^{m}(\phi(b)\ast_{p_{3},p_{2}}^{p_{1}}\phi(a)+x)]\ast_{n,p_{3}}^{m}\phi(u)\\
& &  \ \ \ \ \ - [(\phi(c)\ast_{p_{2},p_{1}}^{m}\phi
(b)+y)\ast_{p_{3},p_{2}}^{m}\phi(a)]\ast_{n,p_{3}}^{m}\phi(u)
\end{eqnarray*}
for some $x\in O'_{p_1,p_{3}}(V)$ and $y\in O'_{m,p_2}(V).$ By
Lemma 2.7, we have
$$(\phi(c)\ast_{p_{3},p_{1}}^{m}O'_{p_{1},p_{3}}(V))\ast_{n,p_{3}}^{m}\phi(u)\subseteq O_{m,n}(V)
$$
and
\begin{eqnarray*}
& &(O'_{m,p_{2}}(V)\ast_{p_{3},p_{2}}^{m}\phi(a))\ast_{n,p_{3}}^{m}\phi(u)\\
& &\equiv
O'_{m,p_{2}}(V)\ast_{n,p_{2}}^{m}(\phi(a)\ast_{n,p_{3}}^{p_{2}}\phi(u))\
({\rm
mod}\ O_{m,n}(V))\\
& & \subseteq O_{m,n}(V).
\end{eqnarray*}
 Hence by Lemma 2.8, we have
\begin{eqnarray*}
& &\ \ \ \
\phi[u{\ast}_{m,p_{3}}^{n}((a\ast_{p_{1},p_{2}}^{p_{3}}b){\ast}_{m,p_{1}}^{p_{3}}c
-
a{\ast}_{m,p_{2}}^{p_{3}}(b{\ast}_{m,p_{1}}^{p_{2}}c))]\\
 & &\equiv
[\phi(c)\ast_{p_{3},p_{1}}^{m}(\phi(b)\ast_{p_{3},p_{2}}^{p_{1}}\phi(a))-(\phi(c)\ast_{p_{2},p_{1}}^{m}\phi
(b))\ast_{p_{3},p_{2}}^{m}\phi(a)]\ast_{n,p_{3}}^{m}\phi(u) \\&
&\equiv 0 \ ({\rm mod}\ O_{m,n}(V)).
\end{eqnarray*}
 Thus $\phi: A_{n,m}(V)\to A_{m,n}(V)$ is a
well defined bimodule isomorphism.
 \qed

It is proved in \cite{DLM3} that the identity map on $V$ induces
epimorphism of associative algebras from $A_n(V)$ to $A_m(V)$, for
$m\leq n.$ A similar result holds here.
\begin{prop} Let $m,n,l$ be nonnegative integers such that $m-l,n-l$ are
nonnegative. Then $A_{n-l,m-l}(V)$ is an
$A_n(V)$-$A_m(V)$-bimodule and the identity map on $V$ induces an
epimorphism of $A_n(V)$-$A_m(V)$-bimodules from $A_{n,m}(V)$ to
$A_{n-l,m-l}(V).$
\end{prop}

\pf It is good enough to prove the results for $l=1.$ We first
show that $u\ast_{p_{1},p_{2}}^{p_{3}}v\equiv
u\ast_{p_{1}-1,p_{2}-1}^{p_{3}-1}v$ {\rm
mod}$O'_{p_{3}-1,p_{1}-1}(V)$, for $p_{1}, p_{2}, p_{3}\in{\mathbb
Z}_{+}$. Let $u$ be homogeneous. Then
\begin{eqnarray*}
&
&u\ast_{p_{1},p_{2}}^{p_{3}}v=\sum_{i=0}^{p_{2}}{p_{1}+p_{3}-p_{2}+i\choose
i}
(-1)^i \Res_zY(u,z)v\frac{(1+z)^{{\wt}\,u+p_{1}-1}}{z^{p_{1}+p_{3}-p_{2}+i}}\\
& & \ \ \ \ \ +\sum_{i=0}^{p_{2}}{p_{1}+p_{3}-p_{2}+i\choose
i}(-1)^i
\Res_zY(u,z)v\frac{(1+z)^{{\wt}\,u+p_{1}-1}}{z^{p_{1}+p_{3}-p_{2}+i+1}}\\
& &\equiv \sum_{i=0}^{p_{2}-1}{p_{1}+p_{3}-p_{2}+i\choose i}(-1)^i
\Res_zY(u,z)v\frac{(1+z)^{{\wt}\,u+p_{1}-1}}{z^{p_{1}+p_{3}-p_{2}+i}}\ \ \\
& &\ \ \ \ \  +\sum_{i=0}^{p_{2}-2}{p_{1}+p_{3}-p_{2}+i\choose
i}(-1)^i\Res_zY(u,z)\frac{(1+z)^{{\wt}\,u+p_{1}-1}}{z^{p_{1}+p_{3}-p_{2}+i+1}}\
\pmod{O'_{p_{3}-1,p_{1}-1}(V)}\\
&
&=\left[\sum_{i=0}^{p_{2}-1}(-1)^{i}{p_{1}+p_{3}-1-p_{2}+i\choose
i}+\sum_{i=1}^{p_{2}-1}(-1)^{i}{p_{1}+p_{3}-1-p_{2}+i\choose
i-1}\right]\\
& & \ \ \ \ \ \ \cdot\Res_zY(u,z)v\frac{(1+z)^{\wt
u+p_{1}-1}}{z^{p_{1}+p_{3}-p_{2}+i}}\\
& &\ \ \ \ \
+\sum_{i=1}^{p_{2}-1}(-1)^{i-1}{p_{1}+p_{3}-1-p_{2}+i\choose
i-1}\Res_zY(u,z)v\frac{(1+z)^{\wt
u+p_{1}-1}}{z^{p_{1}+p_{3}-p_{2}+i}}\\
& &=u\ast_{p_{1}-1,p_{2}-1}^{p_{3}-1}v.
\end{eqnarray*}

In particular, we have $u*^n_mv\equiv u*^{n-1}_{m-1}v$ mod
$O'_{n-1,m-1}(V)$ and $u\bar*^n_mv\equiv u\bar*^{n-1}_{m-1}v$ mod
$O'_{n-1,m-1}(V).$ It remains to prove that $O_{n,m}(V)\subset
O_{n-1,m-1}(V).$ Clearly, $O'_{n,m}(V)\subset O'_{n-1,m-1}(V)$ by
Lemma \ref{l2.4}. Recall from \cite{DLM3} that $O_m(V)\subset
O_{m-1}(V)$ for any $m.$ Using the relation
$u\ast_{p_{1},p_{2}}^{p_{3}}v\equiv
u\ast_{p_{1}-1,p_{2}-1}^{p_{3}-1}v$ mod $O'_{p_{3}-1,p_{1}-1}(V)$
and Lemma 2.7 and the definition of $O_{n,m}(V)$, one can easily
show that $O''_{n,m}(V)\subset O_{n-1,m-1}(V),$
$(V\ast_{p_{1}}^{n}O_{p_{1}}(V))\ast_{m,p_{1}}^{n}V\subset
O_{n-1,m-1}(V).$ \qed
\\

We next study the tensor product
$A_{n,p}(V)\otimes_{A_p(V)}A_{p,m}(V)$ which is an
$A_n(V)$-$A_m(V)$-bimodule.
\begin{prop}\label{p3.4} Define the linear map $\psi$: $A_{n,p}(V)\otimes_{A_p(V)}A_{p,m}(V)\rightarrow
A_{n,m}(V)$ by
$$
\psi(u\otimes v)=u\ast_{m,p}^{n}v,$$ for $u\otimes v\in
A_{n,p}(V)\otimes_{A_p(V)}A_{p,m}(V)$. Then $\psi$ is an
$A_{n}(V)-A_{m}(V)$- bimodule homomorphism from
$A_{n,p}(V)\otimes_{A_p(V)}A_{p,m}(V)$ to $A_{n,m}(V)$.
\end{prop}

\pf First we prove that $\psi$ is well defined. By Lemma 2.8, if
$u\in O_{n,p}(V)$ or $v\in O_{p,m}(V)$, then $\psi(u\otimes v)=0$.
For $u\in A_{n,p}(V), w\in A_{p}(V), v\in A_{p,m}(V)$, we have
$$
\psi((u\ast_{p,p}^{n}w)\otimes
v)=(u\ast_{p,p}^{n}w)\ast_{m,p}^{n}v$$
$$
=u\ast_{m,p}^{n}(w\ast_{m,p}^{p}v)=\psi(u\otimes(w\ast_{m,p}^{p}v)).$$
Therefore $\psi$ is well defined.

For $a\in A_{n}(V)$, $b\in A_{m}(V)$ and $u\otimes v\in
A_{n,p}(V)\otimes_{A_p(V)}A_{p,m}(V)$, by the definition of
$A_{n,m}(V)$, we have
$$
a\bar{\ast}_{m}^{n}\psi(u\otimes
v)=a\bar{\ast}_{m}^{n}(u\ast_{m,p}^{n}v)=(a\bar{\ast}_{p}^{n}u)\ast_{m,p}^{n}v=\psi((a\bar{\ast}_{p}^{n}u)\otimes
v);$$
$$
\psi(u\otimes
v)\ast_{m}^{n}b=(u\ast_{m,p}^{n}v)\ast_{m}^{n}b=u\ast_{m,p}^{n}(v\ast_{m}^{p}b)=\psi(u\otimes
(v\ast_{m}^{p}b)).$$ \qed

It is worthy to point out that the map $\psi$ is not surjective in
general. For example, if $V=V^{\natural}$ is the moonshine vertex
operator algebra constructed in \cite{FLM} then $V$ is rational
(see \cite{D}, \cite{DGH} and \cite{M}) and $V_1=0.$ Thus by Theorem
\ref{t4.10} below, $A_{2,1}(V)=0$ and $A_{2,2}(V)\ne 0.$ This
shows that $A_{2,1}(V)\otimes_{A_1(V)}A_{1,2}(V)=0$ and $\psi:
A_{2,1}(V)\otimes_{A_1(V)}A_{1,2}(V)\to A_2(V)$ is the zero map.

\section{Connection to representation theory}
\def\theequation{4.\arabic{equation}}
\setcounter{equation}{0}

Let $M=\bigoplus_{n\in{\mathbb Z}_{+}}M(n)$ be an admissible $V$-module
such that $M(0)\ne 0$ (cf. \cite{DLM2}).
 For $u\in V$, define $o_{n,m}(u): M(m)\mapsto M(n)$ by
$$o_{n,m}(u)w=u_{wtu+m-n-1}w,$$
for homogeneous $u\in
V$ and $w\in M(m)$ where $u_{wtu+m-n-1}$ is the component
operator of $Y_M(u,z)=\sum_{n\in\Z}u_nz^{-n-1}.$

The following lemma gives the representation theory reason for
Proposition \ref{p3.4}.

\begin{lem}\label{l4.1} Let $a,b\in V$, $m, n, p\in{\mathbb Z}_{+}$ and $w\in M(m)$, then
$$
o_{n,m}(a\ast_{m,p}^{n}b)w=o_{n,p}(a)o_{p,m}(b)w.$$
\end{lem}

\pf We have the following computation on $M(m):$
\begin{eqnarray*}
& &o_{n,m}(a\ast_{m,p}^{n}b)
 =o_{n,m}\left(\sum\limits_{i=0}^{p}(-1)^{i}{n+m-p+i \choose i}
{\rm Res}_{z}\frac{(1+z)^{wta+m}}
{z^{n+m-p+i+1}}Y(a,z)b\right)\\
 & &=\sum\limits_{i=0}^{p}(-1)^{i}{n+m-p+i \choose
i}\sum\limits_{j=0}^{wta+m}{ wta+m \choose
j}(a_{j-n-m+p-i-1}b)_{wta+wtb-j+i+2m-p-1}
\\
& &=\sum\limits_{i=0}^{p}(-1)^{i}{n+m-p+i \choose
i}\sum\limits_{j=0}^{wta+m}{ wta+m \choose j}\\
& &\ \ \ \cdot {\rm Res}_{z_{2}}{\rm
Res}_{z_{1}-z_{2}}(z_{1}-z_{2})^{j-n-m+p-i-1}z_{2}^{wta+wtb-j+i+2m-p-1}
Y(Y(a,z_{1}-z_{2})b,z_{2})\\
& &=\sum\limits_{i=0}^{p}(-1)^{i}{n+m-p+i\choose i}\\
& &\ \ \ \cdot {\rm
Res}_{z_{2}}{\rm Res}_{z_{1}-z_{2}}z_{1}^{wta+m}(z_{1}-z_{2})^{-n-m+p-i-1}z_{2}^{wtb+i+m-p-1}Y(Y(a,z_{1}-z_{2})b,z_{2})\\
& &=\sum\limits_{i=0}^{p}(-1)^{i}{n+m-p+i \choose
i}\sum\limits_{j=0}^{\infty}{ -n-m+p-i-1\choose j}\\
& &\ \ \ \cdot{\rm Res}_{z_{1}}{\rm Res}_{z_{2}}
z_{1}^{-n-m+p-i-1-j}(-z_{2})^{j}z_{1}^{wta+m}z_{2}^{wtb+i+m-p-1}Y(a,z_{1})Y(b,z_{2})
\\
& &\ \ \ \ \  -\sum\limits_{i=0}^{p}(-1)^{i}{ n+m-p+i\choose
i}\sum\limits_{j=0}^{\infty}{-n-m+p-i-1\choose j}\\
& &\ \ \ \cdot {\rm Res}_{z_{2}}{\rm Res}_{z_{1}}
(-z_{2})^{-n-m+p-i-1-j}z_{1}^{j}z_{1}^{wta+m}z_{2}^{wtb+i+m-p-1}Y(b,z_{2})Y(a,z_{1})
\\
& &=\sum\limits_{i=0}^{p}(-1)^{i}{n+m-p+i \choose
i}\sum\limits_{j=0}^{\infty}{ -n-m+p-i-1\choose j}(-1)^{j}\\
& & \ \ \ \ \ \  \cdot a_{wta-n+p-i-j-1}b_{wtb+i+j+m-p-1}\\
&  &\ \ \ -\sum\limits_{i=0}^{p}(-1)^{i}{n+m-p+i \choose
i}\sum\limits_{j=0}^{\infty}{ -n-m+p-i-1\choose j}\\
&  &\ \ \ \ \ \ \cdot(-1)^{n+m-p+j+i+1}
 b_{wtb-n-2-j}a_{wta+m+j}\\
 & &=\sum\limits_{k=0}^{\infty}(-1)^{k}\sum\limits_{i=0}^{p}{n+m-p+i
\choose i}{ -n-m+p-i-1\choose k-i}\\ & &
\ \ \ \ \ \cdot a_{wta-n+p-k-1}b_{wtb+m-p+k-1}\\
& & \ \ \ \ \  -\sum\limits_{i=0}^{p}(-1)^{i}{n+m-p+i \choose
i}\sum\limits_{j=0}^{\infty}{ -n-m+p-i-1\choose j}\\ & & \ \ \ \ \
\cdot (-1)^{n+m-p+j+i+1}b_{wtb-n-2-j}a_{wta+m+j}.
\end{eqnarray*}
Note that $a_{wta+m+j}=b_{wtb+m-p+k-1}=0$ on $M(m)$ for $j\geq 0$
and $k>p.$ Also,
$$
\sum\limits_{i=0}^{k}{ n+m-p+i\choose i}{ -n-m+p-i-1\choose
k-i}=0, \ k=1,2,\cdots,p.
$$
The proof is complete.
\qed

 \begin{coro}\label{c4.2} \ For $u,v\in V$ and $w\in M(m)$, we have
$$o_{n,m}(u\ast_{m}^{n}v)w=o_{n,m}(u)o_{m,m}(v)w, \ \
o_{n,m}(u\bar{\ast}_{m}^{n}v)w=o_{n,n}(u)o_{n,m}(v)w.$$
\end{coro}

Let $W$ be a weak $V$-module and $m\in{\mathbb Z}_{+}.$ Following \cite{DLM2}
we define $$
\Omega_{m}(W)=\{w\in W|u_{wtu-1+k}w=0, {\rm for \ all \
homogeneous } \ u\in V \ {\rm and} \ k>m\}.$$ By Corollary \ref{c4.2}
or \cite{DLM3} we have:
\begin{theorem}  Let $W$ be a weak $V$-module. Then
$\Omega_{m}(W)$ is an $A_m(V)$-module such that $a+O_m(V)$
acts as $o(a)=a_{wta-1}$ for homogeneous $a\in V$.
\end{theorem}

We also have the following theorem from \cite{DLM3}.
\begin{theorem} Let $M=\bigoplus_{k=0}^{\infty}M(k)$ be an admissible
$V$-module such that $M(0)\ne 0.$ Then

(1) $\bigoplus_{k=0}^{m}M(k)$ is an $A_m(V)$-submodule of $\Omega_m(M).$ Furthermore, if $M$ is irreducible,
 $\Omega_m(M)=\bigoplus_{k=0}^{m}M(k).$

(2) Each $M(k)$ is an $A_m(V)$-submodule for $k=0,...,m.$ Moreover, if $M$ is
irreducible, each $M(k)$ is a simple $A_m(V)$-module.
\end{theorem}

We are now in a position to understand the representation theory meaning
of $A_{n,m}(V).$ First observe that
${\rm Hom}_{\C}(M(m),M(n))$ is an $A_n(V)$-$A_m(V)$-bimodule such that
$(a\cdot f\cdot b)(w)=af(bw)$ for $a\in A_n(V), b\in A_m(V),$ $f\in  {\rm Hom}_{\C}(M(m),M(n))$ and $w\in M(m).$
Set $o_{n,m}(V)=\{o_{n,m}(v)|v\in V\}.$
\begin{prop} The $o_{n,m}(V)$ is an $A_n(V)$-$A_m(V)$-subbimodule of
${\rm Hom}_{\C}(M(m),M(n))$ and $v\mapsto o_{n,m}(v)$ for $v\in V$ induces
an $A_n(V)$-$A_m(V)$-bimodule epimorphism from $A_{n,m}(V)$ to $o_{n,m}(V).$
\end{prop}

\pf Clearly, $o_{n,m}(V)$ is an $A_n(V)$-$A_m(V)$-subbimodule of
${\rm Hom}_{\C}(M(m),M(n))$ by Corollary \ref{c4.2}. The same
corollary also shows that the map $v\mapsto o_{n,m}(v)$ is an
$A_n(V)$-$A_m(V)$-bimodule epimorphism if we can prove that
$o_{n,m}(c)=0$ for $c\in O_{n,m}(V).$ First let  $c\in
O'_{n,m}(V).$ If $c=L(-1)u+(L(0)+m-n)u,$  $o_{n,m}(c)=0$ is clear.
If $c=u\circ^n_mv,$ then by Lemma \ref{l2.2} we can assume that
$c=a\bar *_m^nb$ for some $a\in O_n(V)$ and $b\in V.$ Since
$o(a)=0$ we see from Corollary \ref{c4.2} that $o_{n,m}(c)=0.$
Next we assume that $c\in O''_{n,m}(V).$ Using Lemma \ref{l4.1}
repeatedly shows that $o_{n,m}(c)=0.$ Finally take $c\in
O'''_{n,m}(V).$ Again by Lemma \ref{l4.1} and Corollary
\ref{c4.2}, $o_{n,m}(c)=0,$ as desired. \qed

It is proved in \cite{DLM3} that for any given $A_m(V)$-module $U$
which cannot factor through $A_{m-1}(V)$ there is a unique
admissible $V$-module $\bar M(U)=\bigoplus_{k=0}^{\infty}\bar
M(U)(k)$ of Verma type such that $\bar M(U)(m)=U.$ The
construction of  $\bar M(U)$ is implicit in \cite{DLM3}.  We now
recover this result and give an explicit construction of  $\bar
M(U)$ by using the bimodules $A_{n,m}(V).$ As a byproduct of this
construction we can determine the structure of $A_{n,m}(V)$
explicitly if $V$ is rational. We also expect that this new
construction of  $\bar M(U)$ will help our further study of
representation theory of $V.$

We need the following result on the relation
between  $A_{n,n}(V)$ and $A_n(V):$

\begin{prop}\label{pa4.5} For any $n\geq 0,$ the $A_n(V)$ and $A_{n,n}(V)$ are the same.
\end{prop}

\pf It is good enough to prove that any $A_n(V)$-module $U$ is
also an $A_{n,n}(V)$-module. Recall from the definition of
$A_n(V)$ and $A_{n,n}(V)$ that $A_n(V)=V/O'_{n,n}(V)$ and
$A_{n,n}(V)=V/O_{n,n}(V)$ where
$O_{n,n}(V)=O'_{n,n}(V)+O''_{n,n}(V)+O'''_{n,n}(V).$ So we have to
prove that $O''_{n,n}(V)+O'''_{n,n}(V)$ acts on $U$ trivially. By
Theorem 4.1 of \cite{DLM3}, there exists an admissible $V$-module
$M=\bigoplus_{k=0}^{\infty}M(k)$ such that $M(n)=U.$ We can not
assume and do not need to assume that $M(0)\ne 0.$ Note that the
action of $A_n(V)$ on $U\subset M $ comes from the $A_n(V)$-module
structure. This is, for $v\in A_n(V),$ $o(v)=o_{n,n}(v)$ is the
module action of $A_n(V)$ on $U.$ By Lemma \ref{l4.1} we
immediately see that $O''_{n,n}(V)+O'''_{n,n}(V)=0$ on $U,$ as
desired.  \qed

Let $U$ be an $A_{m}(V)$-module which can not factor through
$A_{m-1}(V)$. Set
$$
M(U)=\bigoplus_{n\in{\mathbb Z}_{+}}A_{n,m}(V)\otimes_{A_{m}(V)}
U.$$
 Then $M(U)$
is naturally ${\mathbb Z}_{+}$-graded such that
$M(U)(n)=A_{n,m}(V)\otimes_{A_{m}(V)}U.$ By Proposition
\ref{pa4.5}, $M(U)(m)$ and $U$ are isomorphic $A_{m}(V)$-modules.
The $M(U)$ will be proved to be the $\bar M(U)$ defined in
\cite{DLM3}.

Recall Proposition \ref{p3.4}. For $u\in V,$ and $p,n\in {\mathbb
Z}$, define an operator $u_{p}$ from $M(U)(n)$ to
$M(U)(n-wtu-p-1)$ by
$$
u_{p}(v\otimes
w)=\left\{\begin{array}{rl}(u\ast_{m,n}^{wtu-p-1+n}v)\otimes w,
\quad {\rm if} \ wtu-1-p+n\geq 0, \\  0, \quad {\rm if} \
wtu-1-p+n< 0,
\end{array}\right.
$$
for $v\in A_{n,m}(V)$ and $w\in U$.
We need to prove that $u_{p}$ is well defined. Let $v\in
O_{n,m}(V)$ and $w\in U.$ By Lemma 2.8, $u\ast_{m,n}^{wtu-p-1+n}v\in
V\ast_{m,n}^{wtu-p-1+n}O_{n,m}(V)\subseteq O_{wtu-p-1+n,m}(V)$, so
we have $u_{p}(v\otimes w)=0$. Now let $a\in A_{m}(V)$, $v\in
A_{n,m}(V)$, $w\in U$. Then
\begin{eqnarray*}
& &u_{p}((v\ast_{m}^{n}a)\otimes
w)=(u\ast_{m,n}^{wtu-p-1+n}(v\ast_{m}^{n}a))\otimes w\\
& &\ \ \ \ =((u\ast_{m,n}^{wtu-p-1+n}v)\ast_{m}^{wtu-p-1+n}a)\otimes
w\\
& &\ \ \ \ =(u\ast_{m,n}^{wtu-p-1+n}v)\otimes a\cdot w=u_{p}(v\otimes
a\cdot w).
\end{eqnarray*}
Thus $u_{p}$ is well defined. Set
$$Y_{M}(u,z)=\sum\limits_{p\in{\mathbb Z}}u_{p}z^{-p-1}.$$

\begin{lem}\label{l4.6} For homogeneous $u\in V$ , $v\otimes w\in A_{n,m}(V)\otimes_{A_{m}(V)} U$ and
$p\in{\mathbb Z}$, we have

(1) $u_{p}(v\otimes w)=0$, for $p$ sufficiently large;

(2) $Y_{M}(\1,z)={\rm id}.$
\end{lem}

\pf (1) is clear and we only need to deal with (2). By the
definition of $u_{p}$, we have
\begin{eqnarray*}
& &\1_{p}(v\otimes
w)=\sum\limits_{i=0}^{n}(-1)^{i}{-p-1+m+i\choose i}{\rm
Res}_{z}\frac{(1+z)^{m}}{z^{-p+m+i}}(Y(u,z)v)\otimes
w\\
& &=\sum\limits_{i=0}^{n}(-1)^{i}{-p-1+m+i\choose
i}\sum\limits_{j=0}^{m}{m\choose j}(\1_{p-i-j}v)\otimes w
\end{eqnarray*}
Thus $\1_{p}(v\otimes w)=0$ if $p<-1$ and $\1_{-1}(v\otimes w)=v\otimes w.$
Clearly, $\1_{p}(v\otimes w)=0$ if $p\geq n$. If $-1<p<n$, then
\begin{eqnarray*}
 &
&\1_{p}(v\otimes w)=\sum\limits_{i=0}^{p+1}(-1)^{i}{m-p-1+i\choose
i}{m\choose
p+1-i}v\otimes w\\
& &=\sum\limits_{i=0}^{p+1}(-1)^{i+p+1}{m\choose p+1}{p+1\choose
i}v\otimes w\\
& &=0.
\end{eqnarray*}
That is,  $Y_{M}(\1,z)={\rm id}.$ \qed

We next have the commutator formula:
\begin{lem}\label{l4.7} For $a,b\in V$, we have
$$
[Y_{M}(a,z_{1}),Y_{M}(b,z_{2})]={\rm
Res}_{z_{0}}z_{2}^{-1}\delta(\frac{z_{1}-z_{0}}{z_{2}})Y_{M}(Y(a,z_{0})b,z_{2}),$$
or equivalently, for $p,q\in\Z$
$$[a_p,b_q]=\sum_{i\geq 0}{p\choose i}(a_ib)_{p+q-i}.$$
\end{lem}

\pf Recall Lemma \ref{l2.3} and the definition of $A_{n,m}(V).$
For $p,q\in{\mathbb Z}$ and $v\otimes w\in
A_{n,m}(V)\otimes_{A_{m}(V)} U$  we need to prove that
$$ a_{p}b_{q}(v\otimes w)-b_{q}a_{p}(v\otimes w)=\sum\limits_{i=0}^{\infty}{p\choose i}(a_ib)_{p+q-i}(v\otimes
w).$$ It $wt a+wt b-p-q-2+n<0,$ this is clear from the definition
of the actions. So we now assume that $wt a+wt b-p-q-2+n\geq 0.$
If $wta-p-1+n, wtb-q-1+n\geq 0,$ then by Lemma \ref{l2.3} we have
\begin{eqnarray*}
& &\ \ \ \  a_{p}b_{q}(v\otimes w)-b_{q}a_{p}(v\otimes w)\\
& &=a_{p}(b\ast_{m,n}^{wtb-q-1+n}v)\otimes
w-b_{q}(a\ast_{m,n}^{wta-p-1+n}v)\otimes w\\
&
&=\left(a\ast_{m,wtb-q-1+n}^{wta+wtb-p-q-2+n}(b\ast_{m,n}^{wtb-q-1+n}v)\right)\otimes
w\\
& & \ \ \ \ \   -\left(b\ast_{m,wta-p-1+n}^{wta+wtb-p-q-2+n}(a\ast_{m,n}^{wta-p-1+n}v)\right)\otimes w\\
&
&=\left((a\ast_{n,wtb-q-1+n}^{wta+wtb-p-q-2+n}b)\ast_{m,n}^{wta+wtb-p-q-2+n}v\right)\otimes
w\\
& & \ \ \ \ \  -\left((b\ast_{n,wta-p-1+n}^{wta+wtb-p-q-2+n}a)\ast_{m,n}^{wta+w
b-p-q-2+n}v\right)\otimes w\\
& &=\left(({\rm
Res}_{z}(1+z)^{p}Y_{M}(a,z)b)\ast_{m,n}^{wta+wtb-p-q-2+n}v\right)\otimes
w\\
& &=(\sum\limits_{i=0}^{\infty}{p\choose
i}\left(a_ib)\ast_{m,n}^{wta+wtb-p-q-2+n}v\right)\otimes
w\\
& &=\sum\limits_{i=0}^{\infty}{p\choose i}(a_ib)_{p+q-i}(v\otimes
w).
\end{eqnarray*}

It remains to prove the result for $wta-p-1+n\geq 0, wtb-q-1+n<0$
or $wta-p-1+n<0, wtb-q-1+n\geq 0.$ If $wta-p-1+n<0, wtb-q-1+n\geq
0$ then $b_{q}a_{p}(v\otimes w)=0$ and
\begin{eqnarray*}
& &a_{p}b_{q}(v\otimes w)-b_{q}a_{p}(v\otimes w)=\left((a\ast_{n,wtb-q-1+n}^{wta+wtb-p-q-2+n}b)
\ast_{m,n}^{wtp+wtb-p-q-2+n}v\right)\otimes
w\\
& &\ \ \ \ =\left(({\rm
Res}_{z}(1+z)^{p}Y_{M}(a,z)b)\ast_{m,n}^{wta+wtb-p-q-2+n}v\right)\otimes
w\\
& &\ \ \ \ =\sum\limits_{i=0}^{\infty}{p\choose i}(a_ib)_{p+q-i}(v\otimes
w)
\end{eqnarray*}
where we have used Lemma \ref{la4.8} below. Similarly, the result
holds for $wta-p-1+n\geq 0, wtb-q-1+n<0.$ \qed

\begin{lem}\label{la4.8}
Let $u,v\in V$ and $m\in{\mathbb Z}_{+}$, $p_1,p_2\in\Z$ such that
$p_1+p_2-m\geq 0.$ If $p_1\geq 0, p_2<0$, then
$$u{\ast}_{m,p_{1}}^{p_{1}+p_{2}-m}v-{\rm Res}_{z}(1+z)^{wtu-1+m-p_{2}}Y(u,z)v\in
O'_{p_{1}+p_{2}-m,m}(V)$$ and if $p_1<0, p_2\geq 0$, then
$$-v\ast_{m,p_{2}}^{p_{1}+p_{2}-m}u-{\rm Res}_{z}(1+z)^{wtu-1+m-p_{2}}Y(u,z)v\in
O'_{p_{1}+p_{2}-m,m}(V).$$
\end{lem}

\pf This lemma is similar to Lemma \ref{l2.3} where both $p_1$ and
$p_2$ are nonnegative. If $p_2<0$ or $p_1<0$, then
$v\ast_{m,p_{2}}^{p_{1}+p_{2}-m}u$ or $
u{\ast}_{m,p_{1}}^{p_{1}+p_{2}-m}v$ is not defined. But we do need
a version of Lemma \ref{l2.3} with either $p_1<0$ or $p_2<0$ in
the proof of previous lemma.

First we assume that  $p_1\geq 0, p_2<0.$ Then $-p_2-1\geq 0.$ From the
definition, we have
$$u{\ast}_{m,p_{1}}^{p_{1}+p_{2}-m}v=\sum_{i=0}^{p_1}{-p_2-1\choose i}
\Res_zY(u,z)v\frac{(1+z)^{wtu+m}}{z^{p_2+i+1}}.$$ Since
$\Res_zY(u,z)v\frac{(1+z)^{wtu+m}}{z^{p_2+i+1}}\in
O'_{p_1+p_2-m,m}(V)$ if $i>p_1$  we see that
\begin{eqnarray*}
& &u{\ast}_{m,p_{1}}^{p_{1}+p_{2}-m}v\equiv
\sum_{i=0}^{-p_2-1}{-p_2-1\choose i}
\Res_zY(u,z)v\frac{(1+z)^{wtu+m}}{z^{p_2+i+1}}.\\
& &\ \ \ \ =\Res_zY(u,z)v\frac{(1+z)^{wtu+m}}{z^{p_2+1}}
(1+\frac{1}{z})^{-p_2-1}\\
& &\ \ \ \ = \Res_zY(u,z)v(1+z)^{wtu+m-p_2-1}.
\end{eqnarray*}
So in this case we have done.

If  $p_1<0, p_2\geq 0$ then the result in the first case
gives
$$v\ast_{m,p_{2}}^{p_{1}+p_{2}-m}u\equiv \Res_zY(v,z)u(1+z)^{wtv+m-p_1-1}$$
modulo $O'_{p_{1}+p_{2}-m,m}(V).$
Using the identity
$$Y(v,z)u\equiv(1+z)^{-wtu-wtv-2m+p_{1}+p_{2}}Y(u,\frac{-z}{1+z})v$$
modulo $O'_{p_{1}+p_{2}-m,m}(V)$ we see that
\begin{eqnarray*}
& & \Res_zY(v,z)u(1+z)^{wtv+m-p_1-1}\equiv
\Res_zY(u,\frac{-z}{1+z})v(1+z)^{-wtu-m+p_{2}-1}\\
& &\ \ \ \ =-\Res_zY(u,z)v(1+z)^{wtu+m-p_{2}-1}.
\end{eqnarray*}
The proof is complete.
\qed

\begin{lem}\label{l4.8} For $i\in{\mathbb Z}_{+},$ we have
\begin{eqnarray*}
{\rm
Res}_{z_{0}}z_{0}^{i}(z_{0}+z_{2})^{wta+n}Y_{M}(a,z_{0}+z_{2})Y_{M}(b,z_{2})
={\rm
Res}_{z_{0}}z_{0}^{i}(z_{2}+z_{0})^{wta+n}Y_{M}(Y(a,z_{0})b,z_{2})
\end{eqnarray*}
on $M(U)(n)$ for $n\in\Z_+.$
\end{lem}

\pf Since $a_{wta+n}=0$ on $M(U)(n),$ we have
$$
{\rm
Res}_{z_{1}}(z_{1}-z_{2})^{i}z_{1}^{wta+n}Y_{M}(b,z_{2})Y_{M}(a,z_{1})(v\otimes
w)=0.$$ Thus on $M(U)(n),$ we have
\begin{eqnarray*}
& &{\rm
Res}_{z_{0}}z_{0}^{i}(z_{0}+z_{2})^{wta+n}Y_{M}(a,z_{0}+z_{2})Y_{M}(b,z_{2})\\
& &={\rm
Res}_{z_{1}}(z_{1}-z_{2})^{i}z_{1}^{wta+n}(Y_{M}(a,z_{1})Y_{M}(b,z_{2})-Y_{M}(b,z_{2})Y_{M}(a,z_{1}))\\
& &={\rm
Res}_{z_{1}}(z_{1}-z_{2})^{i}z_{1}^{wta+n}[Y_{M}(a,z_{1}),Y_{M}(b,z_{2})]\\
& &={\rm Res}_{z_{0}}{\rm
Res}_{z_{1}}(z_{1}-z_{2})^{i}z_{1}^{wta+n}z_{2}^{-1}\delta(\frac{z_{1}-z_{0}}{z_{2}})Y_{M}(Y(a,z_{0})b,z_{2})\\
& &={\rm Res}_{z_{0}}{\rm
Res}_{z_{1}}z_{0}^{i}(z_{2}+z_{0})^{wta+n}z_{1}^{-1}\delta(\frac{z_{2}+z_{0}}{z_{1}})Y_{M}(Y(a,z_{0})b,z_{2})\\
& &={\rm
Res}_{z_{0}}z_{0}^{i}(z_{2}+z_{0})^{wta+n}Y_{M}(Y(a,z_{0})b,z_{2}),
\end{eqnarray*}
where we have used Lemma \ref{l4.7}.\qed

\begin{lem}\label{l4.9} For $l\in{\mathbb Z}_{+},$ we have
\begin{eqnarray*}
& &{\rm
Res}_{z_{0}}z_{0}^{-l}(z_{2}+z_{0})^{wta+n}z_{2}^{wtb-n}Y_{M}(Y(a,z_{0})b,z_{2})\\
& &={\rm
Res}_{z_{0}}z_{0}^{-l}(z_{0}+z_{2})^{wta+n}z_{2}^{wtb-n}Y_{M}(a,z_{0}+z_{2})Y_{M}(b,z_{2})
\end{eqnarray*}
on $M(U)(n)$ for $n\geq 0.$
\end{lem}

\pf Take $v\otimes w\in A_{n,m}(V)\otimes_{A_m(V)}U=M(U)(n).$ Then
\begin{eqnarray*} & &{\rm
Res}_{z_{0}}z_{0}^{-l}(z_{2}+z_{0})^{wta+n}z_{2}^{wtb-n}Y_{M}(Y(a,z_{0})b,z_{2})(v\otimes w)\\
& &=\sum_{j\in{\mathbb Z}_{+}}{wta+n\choose
j}z_{2}^{wta+wtb-j}Y_{M}(a_{j-l}b,z_{2})(v\otimes w)\\
& &=\sum_{j\in{\mathbb Z}_{+}}{wta+n\choose
j}\sum\limits_{k\in{\mathbb
Z}_{+}}z_{2}^{-l+k-n+1}(a_{j-l}b)_{wta+wtb-j+l-2-k+n}(v\otimes w)\\
& &=\sum_{k\in{\mathbb
Z}_{+}}z_{2}^{-l+k-n+1}\sum_{j\in{\mathbb
Z}_{+}}{wta+n\choose j}\left((a_{j-l}b)\ast_{m,n}^{k}v\right)\otimes w\\
& &=\sum_{k\in{\mathbb Z}_{+}}z_{2}^{-l+k-n+1}\left(({\rm
Res}_{z}\frac{(1+z)^{wta+n}}{z^{l}}Y(a,z)b)\ast_{m,n}^{k}v\right)\otimes w
\end{eqnarray*}
On the other hand, we have
\begin{eqnarray*}
& &{\rm
Res}_{z_{0}}z_{0}^{-l}(z_{0}+z_{2})^{wta+n}z_{2}^{wtb-n}Y_{M}(a,z_{0}+z_{2})Y_{M}(b,z_{2})(v\otimes
w)\\
& &=\sum\limits_{i\in{\mathbb Z}_{+}}{-l\choose
i}(-1)^{i}a_{wta+n-l-i}z_{2}^{wtb-n+i}Y_{M}(b,z_{2})(v\otimes w)\\
& &=\sum\limits_{i\in{\mathbb Z}_{+}}{-l\choose
i}(-1)^{i}a_{wta+n-l-i}\sum\limits_{j\geq
-n}z_2^{-n+i+j}b_{wtb-1-j}(v\otimes w)\\
& &=\sum_{i\in{\mathbb
Z}_{+}}\sum_{\stackrel{j\geq -n}{l+i+j\geq
1}}{-l\choose
i}(-1)^{i}z_2^{-n+i+j}\left(a\ast_{m,j+n}^{l+i+j-1}(b\ast_{m,n}^{j+n}v)\right)\otimes w\\
& &=\sum\limits_{k\in{\mathbb
Z}_{+}}\sum\limits_{j=-n}^{k+1-l}z_2^{-n+k+1-l}(-1)^{k+1-j-l}{-l\choose
k+1-j-l}\left((a\ast_{n,j+n}^{k}b)\ast_{m,n}^{k}v\right)\otimes w\\
& &=\sum\limits_{k\in{\mathbb
Z}_{+}}\sum\limits_{j=0}^{k+1+n-l}z_2^{-n+k+1-l}(-1)^{k+1+n-j-l}{-l\choose
k+1+n-j-l}\left((a\ast_{n,j}^{k}b)\ast_{m,n}^{k}v\right)\otimes w.
\end{eqnarray*}
So it is enough to prove that
$$\sum\limits_{j=0}^{k+1+n-l}(-1)^{k+1+n-j-l}{-l\choose
k+1+n-j-l}a\ast_{n,j}^{k}b={\rm
Res}_{z}\frac{(1+z)^{wta+n}}{z^{l}}Y(a,z)b.$$
 But
\begin{eqnarray*}
& &\sum\limits_{j=0}^{k+1+n-l}(-1)^{k+1+n-j-l}{-l\choose
k+1+n-j-l}a\ast_{n,j}^{k}b\\
& &=\sum\limits_{j=0}^{k+1+n-l}(-1)^{k+1+n-j-l}{-l\choose
k+1+n-j-l}\sum\limits_{i=0}^{j}(-1)^{i}{k+n-j+i\choose i}\\
& & \ \ \ \ \ \ {\rm
Res}_{z}\frac{(1+z)^{wta+n}}{z^{n+k-j+i+1}}Y(a,z)b\\
& &=\sum\limits_{j=0}^{k+1+n-l}(-1)^{j}{-l\choose
j}\sum\limits_{i=0}^{k+n+1-l-j}(-1)^{i}{l+j+i-1\choose i}{\rm
Res}_{z}\frac{(1+z)^{wta+n}}{z^{l+j+i}}Y(a,z)b.
\end{eqnarray*}
By Proposition \ref{p5.3} in Appendix, we see that
$$
\sum\limits_{j=0}^{k+1+n-l}(-1)^{j}{-l\choose
j}\sum\limits_{i=0}^{k+n+1-l-j}(-1)^{i}{l+j+i-1\choose
i}\frac{1}{z^{l+j+i}}=\frac{1}{z^{l}}.$$
This finishes the proof. \qed

\begin{coro}\label{ca} For $n\in \Z_+,$ we have
\begin{eqnarray*}
(z_{2}+z_{0})^{wta+n}Y_{M}(Y(a,z_{0})b,z_{2})=(z_{0}+z_{2})^{wta+n}Y_{M}(a,z_{0}+z_{2})Y_{M}(b,z_{2})
\end{eqnarray*}
on $M(U)(n).$
\end{coro}

\pf The corollary is an immediate consequence of Lemmas \ref{l4.8}
and \ref{l4.9}. \qed

We now can state another main theorem of this paper.
\begin{theorem}\label{t4.10} Let $U$ be an $A_{m}(V)$-module which can not
factor through $A_{m-1}(V)$, then $M(U)=\bigoplus_{n\in{\mathbb
Z}_{+}}A_{n,m}(V)\otimes_{A_{m}(V)} U$ is an admissible $V$-module
with $M(U)(0)=A_{0,m}(V)\otimes_{A_{m}(V)} U\neq 0$ and with the
following universal property: for any weak $V$-module $W$ and any
$A_{m}(V)$-morphism $\phi: U\rightarrow \Omega_{m}(W)$, there is a
unique morphism $\bar\phi: M(U)\rightarrow W$ of weak $V$-modules
which extends $\phi$.
\end{theorem}

\pf Since $u\bar{\ast}_{m}^{n}1=u$, for $u\in V$ and $1\in
A_{m}(V)$, it follows that $M(U)$ is generated by $U\cong
A_{m,m}(V)\otimes_{A_{m}(V)}U$ by Proposition \ref{pa4.5}. The
fact that $M(U)$ is an admissible $V$-module follows from Lemmas
\ref{l4.6}-\ref{l4.7} and Corollary \ref{ca}. In fact, from the
construction of $M(U)$ and Lemma \ref{l4.6}, every condition in
the definition of an admissible module (cf. \cite{DLM2}) except
the Jacobi identity is satisfied. It is well known that the Jacobi
identity is equivalent to the commutator formula obtained in Lemma
\ref{l4.6} and associativity in Corollary \ref{ca}.  If
$M(U)(0)=A_{0,m}(V)\otimes_{A_{m}(V)} U=0$, then $U$ will be an
$A_{m-1}(V)$-module, a contradiction. So $M(U)(0)\neq 0$.

It remains to prove the universal property of $M(U)$. Define the
linear map: $\bar{\phi}: M(U)\rightarrow W$ by
$$
\bar{\phi}(u\otimes w)=o_{n,m}(u)\phi(w),$$ for $u\in A_{n,m}(V)$
and $w\in U$. To be sure that $\bar\phi$ is well defined, we need
to prove that
$$
\bar{\phi}((u\ast_{m}^{n}v)\otimes w)=\bar{\phi}(u\otimes v\cdot
w),$$ for $u\in A_{n,m}(V)$, $v\in A_{m}(V)$ and $w\in U.$ Indeed,
by  Corollary \ref{c4.2} and the fact that $\phi$ is an
$A_{m}(V)$-morphism, we have
\begin{eqnarray*}
& &\bar{\phi}((u\ast_{m}^{n}v)\otimes
w)=o_{n,m}(u\ast_{m}^{n}v)\phi(w)
=o_{n,m}(u)o(v)\phi(w)\\
& &=o_{n,m}(u)\phi(o(v)w)=\bar{\phi}(u\otimes o(v)w).
\end{eqnarray*}
Hence $\bar{\phi}$ is well defined.

For homogeneous $u\in V$,
$v\otimes w\in A_{n,m}(V)\otimes_{A_{m}(V)}U$ and $p\in{\mathbb
Z}$, by Lemma 4.1, we have
\begin{eqnarray*}
& &\bar{\phi}(u_{p}(v\otimes
w))=\bar{\phi}((u\ast_{m,n}^{wtu-p-1+n}v)\otimes w)\\
&
&=o_{wtu-p-1+n,m}(u\ast_{m,n}^{wtu-p-1+n}v)\phi(w)\\
& &=o_{wtu-p-1+n,n}(u)o_{n,m}(v)\phi(w)
\\
& &=o_{wtu-p-1+n,n}(u)\bar{\phi}(v\otimes
w)\\
& &=u_{p}\bar{\phi}(v\otimes w).
\end{eqnarray*}
This means that $\bar{\phi}$ is a morphism of weak $V$-modules. It
is clear that $\bar{\phi}$ extends $\phi$. \qed

From the universal property of $M(U)$ we immediately have
\begin{coro}  The  $M(U)$ is isomorphic to the
admissible $V$-module $\bar{M}_{m}(U)$ constructed in \cite{DLM3}.
\end{coro}

\begin{remark} In the case $U=A_m(V),$ then $M(U)=\bigoplus_{n\geq 0}A_{n,m}(V)$
is an admissible $V$-module for any $m\geq 0.$ On the other hand, it is easy
to see that $\bigoplus_{n\geq 0}A_n(V)$ is {\em not} an admissible $V$-module.
We certainly expect that the admissible module  $\bigoplus_{n\geq 0}A_{n,m}(V)$
will play a significant role in our further study of representation theory for
vertex operator algebras.
\end{remark}

Finally we study the $A_n(V)$-$A_m(V)$-bimodule structure of
$A_{n,m}(V)$ if $V$ is rational.  The  result below is not
surprising from the representation theory point of view. Recall
from \cite{DLM3} that if $V$ is rational then there are only
finitely irreducible admissible $V$-modules up to isomorphisms and
each irreducible admissible module is ordinary.

\begin{theorem} If $V$ is a rational vertex operator
algebra and $W^{j}\!=\!\bigoplus_{n\geq 0}\!W^j(n)$ with
$W^j(0)\neq 0$ for $j=1,2,\cdots, s$ are all the inequivalent
irreducible modules of $V,$ then
$$A_{n,m}(V)\cong \bigoplus_{l=0}^{{\rm
min}\{m,n\}}\left(\bigoplus_{i=1}^{s}{\rm Hom}_{\mathbb
C}(W^i(m-l),W^i(n-l))\right).
$$
\end{theorem}

\pf  Since $V$ is rational, $A_n(V)$ is isomorphic to the direct sum of
full matrix algebras $\bigoplus_{i=1}^s\bigoplus_{k=0}^n{\rm End}_{\C}(W^i(k))$
by Theorem 4.10 of \cite{DLM3}. So as an $A_n(V)$-$A_m(V)$-bimodule,
$$A_{n,m}(V)=\bigoplus_{i,j=1}^s\bigoplus_{0\leq p\leq m, 0\leq q\leq n}
c_{i,j,p,q}{\rm Hom}_{\mathbb
C}(W^i(p),W^j(q))$$
for some nonnegative integers $c_{i,j,p,q}.$ So we need to prove that
$c_{i,j,p,q}=0$ if $i\ne j$ or $(p,q)\ne (m-l,n-l)$ for some $0\leq l\leq {\rm min}\{m,n\}$ and $c_{i,i,m-l,n-l}=1.$

We need a general result. Let $U$ be an irreducible
$A_p(V)$-module which is not an $A_{p-1}(V)$-module. Then
$M(U)=\bigoplus_{k=0}^{\infty}A_{k,p}(V)\otimes_{A_p(V)}U$ is an
admissible $V$-module generated by $U$ by Theorem \ref{t4.10}.
Since $V$ is rational, $M(U)$ is a direct sum of irreducible
$V$-modules. Note that $M(U)(p)=U$ generates an irreducible
submodule of $M(U)$. This shows that $M(U)$ is irreducible.

Consider irreducible admissible $V$-module
$M=\bigoplus_{k=0}^{\infty}A_{k,m}(V)\otimes_{A_m(V)}W^i(m).$ Then
$M$ is isomorphic to $W^i.$ Thus each $M(k)$  is isomorphic to
$W(k)$ as $A_k(V)$-module. In particular,
$A_{n,m}(V)\otimes_{A_m(V)}W^i(m)$ is isomorphic to $W^i(n).$ This
shows that $c_{i,i,m,n}=1$ and all other $c_{i,j,m,q}=0$ if either
$j\ne i$ or $q\ne n.$

Next we consider
$M=\bigoplus_{k=0}^{\infty}A_{k,m}(V)\otimes_{A_m(V)}(W^i(m)+W^i(m-1)).$
Then $M$ is isomorphic to $W^i\bigoplus W^i,$
$M(k)=W^i(k)\bigoplus W^i(k-1)$ for $k>0,$ and $M(0)=W^i(0).$ This
implies that $c_{i,i,m-1,n-1}=1$ and $c_{i,j,m-1,q}=0$ if $j\ne i$
or $q\ne n-1.$ Continuing in this way gives the result. \qed

\section{Appendix}
\def\theequation{5.\arabic{equation}}
\setcounter{equation}{0}

In this section we present several identities involving formal variables
used in the previous sections.

For $m,n\in{\mathbb Z}_{+}$, define
$$
A_{m,n}(z)=\sum\limits_{i=0}^{m}(-1)^{i}{ n+i\choose
i}\frac{(1+z)^{n+1}}{z^{n+i+1}}-\sum\limits_{i=0}^{n}(-1)^{m}{
m+i\choose i}\frac{(1+z)^i}{z^{m+i+1}}.$$
\begin{prop}\label{p5.1} For all $m,n\in{\mathbb Z}_{+}$, $A_{m,n}(z)=1$.
\end{prop}
\pf  We will prove that $A_{m,n}(z)=A_{m,n+1}(z)$.
\begin{eqnarray*}
& &A_{m,n+1}(z) =\sum\limits_{i=0}^{m}(-1)^{i}{ n+i+1\choose
i}\frac{(1+z)^{n+2}}{z^{n+i+2}}-\sum\limits_{i=0}^{n+1}(-1)^{m}{
m+i\choose i}\frac{(1+z)^i}{z^{m+i+1}}\\
& &=\left[\sum\limits_{i=0}^{m}(-1)^{i}{ n+i\choose
i}+\sum\limits_{i=1}^{m}(-1)^{i}{ n+i\choose i-1}\right]
\left[\frac{(1+z)^{n+1}}{z^{n+i+2}}+\frac{(1+z)^{n+1}}{z^{n+i+1}}\right]\\
& &-\sum\limits_{i=0}^{n}(-1)^{m}{ m+i\choose
i}\frac{(1+z)^i}{z^{m+i+1}}-(-1)^{m}{
m+n+1\choose n+1}\frac{(1+z)^{n+1}}{z^{m+n+2}}\\
& &=A_{m,n}(z)+\left[\sum\limits_{i=1}^{m}(-1)^{i}{ n+i\choose
i-1}+\sum\limits_{i=0}^{m-1}(-1)^{i+1}{
n+i+1\choose i}\right]\frac{(1+z)^{n+1}}{z^{n+i+2}}\\
& &-(-1)^{m}{ m+n+1\choose
n+1}\frac{(1+z)^{n+1}}{z^{m+n+2}}+\sum\limits_{i=0}^{m}(-1)^{i}{
n+i\choose i}\frac{(1+z)^{n+1}}{z^{n+i+2}}\\
& &=A_{m,n}(z)+\sum\limits_{i=0}^{m}\left[(-1)^{i+1}{ n+i\choose
i}+(-1)^{i}{n+i\choose i}\right]\frac{(1+z)^{n+1}}{z^{n+i+2}}\\
& &-\left[(-1)^{m}{m+n+1\choose n+1}+(-1)^{m+1}{ m+n+1\choose
n+1}\right]\frac{(1+z)^{n+1}}{z^{m+n+2}}\\\\
& &=A_{m,n}(z).
\end{eqnarray*}
Thus we have $$A_{m,n}(z)=A_{m,0}(z)=\sum_{i=0}^{m}(-1)^{i}\frac{(1+z)}{z^{i+1}}-(-1)^{m}\frac{1}{z^{m+1}}$$
which is clearly equal to 1.\qed

\begin{prop}\label{p5.2} For $m,n\in{\mathbb Z}_{+}$, we have
\begin{equation}
\sum\limits_{i=0}^{n}(-1)^{i}\left(\begin{array}{c}
m+i\\i\end{array}
\right)\left(\sum\limits_{j=0}^{n-i}(-1)^{j}\left(\begin{array}{c}
-m-i-1\\j\end{array}
\right)\sum\limits_{l=0}^{i}\left(\begin{array}{c} i\\l\end{array}
\right)\frac{z_{2}^{j+l}}{z_{1}^{j+i}}-\frac{1}{z_{1}^{i}}\right)=0.
\end{equation}
\end{prop}
\pf Set
$$
D_{n,m}(z_{1},z_{2}) =\sum\limits_{i=0}^{n}(-1)^{i}{ m+i\choose
i}\left(\sum\limits_{j=0}^{n-i}(-1)^{j}{ -m-i-1\choose
j}\sum\limits_{l=0}^{i}{i\choose
l}\frac{z_{2}^{j+l}}{z_{1}^{j+i}}-\frac{1}{z_{1}^{i}}\right). $$
We prove by induction on both $m$ and $n$ that
$D_{n,m}(z_1,z_2)=0.$ It is clear that $D_{0,0}(z_1,z_2)=0.$ We
now assume that $D_{n,m}(z_1,z_2)=0.$ Then
\begin{eqnarray*}
& &D_{n,m+1}(z_{1},z_{2})\!=\!\sum\limits_{i=0}^{n}(-1)^{i}{
m+i+1\choose i}\left(\sum\limits_{j=0}^{n-i}(-1)^{j}{
-m-i-2\choose j}\sum\limits_{l=0}^{i}{ i\choose
l}\frac{z_{2}^{j+l}}{z_{1}^{j+i}}-\frac{1}{z_{1}^{i}}\right)\\
& &=\sum\limits_{i=0}^{n}(-1)^{i}{ m+i\choose
i}\sum\limits_{j=0}^{n-i}(-1)^{j}{ -m-i-1\choose
i}\frac{m+i+j+1}{m+1}\sum\limits_{l=0}^{i}{ i\choose l}\frac{z_{2}^{j+l}}{z_{1}^{j+i}}\\
& &\ \ \ \ \ \ -\sum\limits_{i=0}^{n}(-1)^{i}\left[{ m+i\choose i}
+{
m+i\choose i-1}\right]\frac{1}{z_{1}^{i}}\\
& &=D_{n,m}(z_{1},z_{2})-\sum\limits_{i=0}^{n}(-1)^{i}{
m+i\choose i}\frac{i}{m+1}\frac{1}{z_{1}^{i}}\\
& &\ \ \ \ +\sum\limits_{i=0}^{n}(-1)^{i}{ m+i\choose
i}\sum\limits_{j=0}^{n-i}(-1)^{j}{ -m-i-1\choose j}\frac{i+j}{m+1}\sum\limits_{l=0}^{i}{ i\choose
l}\frac{z_{2}^{j+l}}{z_{1}^{j+i}}.
\end{eqnarray*}
By induction assumption, we have
\begin{eqnarray*}
& &\sum\limits_{i=0}^{n}(-1)^{i}{ m+i\choose
i}\sum\limits_{j=0}^{n-i}(-1)^{j}{ -m-i-1\choose j}\frac{i+j}{m+1}\sum\limits_{l=0}^{i}{ i\choose
l}\frac{z_{2}^{j+l}}{z_{1}^{j+i}}\\
& &\ \ \ \ \ -\sum\limits_{i=0}^{n}(-1)^{i}{
m+i\choose i}\frac{i}{m+1}\frac{1}{z_{1}^{i}}\\
& &=0.
\end{eqnarray*}
Thus $D_{n,m+1}(z_{1},z_{2})=0.$ We now deal with $D_{n+1,m}(z_1,z_2).$
We have the computation:
\begin{eqnarray*}
& &D_{n+1,m}(z_{1},z_{2})\\
& &=\sum\limits_{i=0}^{n+1}(-1)^{i}{ m+i\choose
i}\left(\sum\limits_{j=0}^{n+1-i}(-1)^{j}{ -m-i-1\choose
j}\sum\limits_{l=0}^{i}{i\choose l}\frac{z_{2}^{j+l}}{z_{1}^{j+i}}-\frac{1}{z_{1}^{i}}\right)\\
& &=\sum\limits_{i=0}^{n}(-1)^{i}{ m+i\choose
i}\sum\limits_{j=0}^{n+1-i}(-1)^{j}{ -m-i-1\choose
j}\sum\limits_{l=0}^{i}{i\choose l}\frac{z_{2}^{j+l}}{z_{1}^{j+i}}
\\
& & \ \ \ \ \ +(-1)^{n+1}{ m+n+1\choose
n+1}\sum\limits_{l=0}^{n+1}{
n+1\choose l}\frac{z_{2}^{l}}{z_{1}^{n+1}}\\
& &\ \ \ \ \ -\sum\limits_{i=0}^{n}(-1)^{i}{ m+i\choose
i}\frac{1}{z_{1}^{i}}-(-1)^{n+1}{
m+n+1\choose n+1}\frac{1}{z_{1}^{n+1}}\\
& &=\sum\limits_{i=0}^{n}(-1)^{i}{ m+i\choose i}{-m-i-1\choose
n+1-i}(-1)^{n+1-i}\sum\limits_{l=0}^{i}{
i\choose l}\frac{z_{2}^{n+1-i+l}}{z_{1}^{n+1}}\\
& &\ \ \ \ \ +(-1)^{n+1}{m+n+1\choose n+1}\sum\limits_{l=0}^{n+1}{
n+1\choose
l}\frac{z_{2}^{l}}{z_{1}^{n+1}}\\
& & \ \ \ \ \ -(-1)^{n+1}{ m+n+1\choose n+1}\frac{1}{z_{1}^{n+1}}+D_{n,m}(z_{1},z_{2})\\
& &=\sum\limits_{l'=1}^{n+1}\sum\limits_{i=0}^{n}(-1)^{n+1}{
m+i\choose i}{ -m-i-1\choose n+1-i}{
i\choose l'-n-1+i}\frac{z_{2}^{l'}}{z_{1}^{n+1}}\\
& & \ \ \ \ \ +\sum\limits_{l'=1}^{n+1}(-1)^{n+1}{ m+n+1\choose
n+1}{
n+1\choose l'}\frac{z_{2}^{l'}}{z_{1}^{n+1}}\\
& &=\sum\limits_{l'=1}^{n+1}\sum\limits_{i=0}^{n}(-1)^{i}{
m+n+1\choose n+1}{n+1\choose l'}{
l'\choose n+1-i}\frac{z_{2}^{l'}}{z_{1}^{n+1}}\\
& &\ \ \ \ \ \ +\sum\limits_{l'=1}^{n+1}(-1)^{n+1}{ m+n+1\choose
n+1}{ n+1\choose l'}\frac{z_{2}^{l'}}{z_{1}^{n+1}}\\
& &=0.
\end{eqnarray*}
The Proposition is proved. \qed
\begin{prop}\label{p5.3} Let $k,l\in{\mathbb Z}_{+}$ be such that $k+1-l\geq
0$ and $l\geq 1$, then
$$\sum\limits_{j=0}^{k+1-l}(-1)^{j}{-l \choose
j}\sum\limits_{i=0}^{k+1-l-j}(-1)^{i}{l+i+j-1 \choose
i}\frac{1}{z^{i+j+l}}=\frac{1}{z^{l}}.$$
\end{prop}
\pf We prove the result by induction on $l\geq 1$. For $l=1$, the identity
becomes
$$
\sum\limits_{j=0}^{k}\sum\limits_{i=0}^{k-j}(-1)^{i}{j+i\choose
i}\frac{1}{z^{j+i+1}}=\frac{1}{z}.
$$
It is easy to deduce the above formula by induction on $k\geq 0$.
Now, suppose the identity  is true for $l$. Then
$$
\frac{d}{dz}\left(\sum\limits_{j=0}^{k+1-l}(-1)^{j}{-l \choose
j}\sum\limits_{i=0}^{k+1-l-j}(-1)^{i}{l+i+j-1 \choose
i}\frac{1}{z^{i+j+l}}\right)=\frac{d}{dz}\left(\frac{1}{z^{l}}\right),$$
whose left side is
\begin{eqnarray*}
& &\sum\limits_{j=0}^{k+1-l}(-1)^{j}{-l \choose
j}\sum\limits_{i=0}^{k+1-l-j}(-1)^{i}{l+i+j-1 \choose
i}(-l-j-i)\frac{1}{z^{i+j+l+1}}\\
& &=(-l)\sum\limits_{j=0}^{k+1-l}(-1)^{j}{-l-1 \choose
j}\sum\limits_{i=0}^{k+1-l-j}(-1)^{i}{l+i+j\choose
i}\frac{1}{z^{i+j+l+1}}\\
& &=(-l)\sum\limits_{j=0}^{k-l}(-1)^{j}{-l-1 \choose
j}\sum\limits_{i=0}^{k+1-l-j}(-1)^{i}{l+i+j\choose
i}\frac{1}{z^{i+j+l+1}}\\
& & \ \ \ \ \ -l{-l-1\choose k+1-l}(-1)^{k+1-l}\frac{1}{z^{k+2}}\\
& &=(-l)\sum\limits_{j=0}^{k-l}(-1)^{j}{-l-1 \choose
j}\sum\limits_{i=0}^{k-l-j}(-1)^{i}{l+i+j\choose
i}\frac{1}{z^{i+j+l+1}}\\
& &\ \ \ \ \ -l\sum\limits_{j=0}^{k-l}{-l-1 \choose
j}(-1)^{k+1-l}{k+1\choose k+1-l-j}\frac{1}{z^{k+2}}\\
& & \ \ \ \ \ -l{-l-1\choose k+1-l}(-1)^{k+1-l}\frac{1}{z^{k+2}}\\
& &=(-l)\sum\limits_{j=0}^{k-l}(-1)^{j}{-l-1 \choose
j}\sum\limits_{i=0}^{k-l-j}(-1)^{i}{l+i+j\choose
i}\frac{1}{z^{i+j+l+1}}\\
& &\ \ \ \ \ -l{k+1\choose
l}(-1)^{k+1-l}\sum\limits_{j=0}^{k-l}{k+1-l\choose
j}(-1)^{j}\frac{1}{z^{k+2}}\\
& & \ \ \ \ \ -l{k+1\choose l}\frac{1}{z^{k+2}}\\
& &=(-l)\sum\limits_{j=0}^{k-l}(-1)^{j}{-l-1 \choose
j}\sum\limits_{i=0}^{k-l-j}(-1)^{i}{l+i+j\choose
i}\frac{1}{z^{i+j+l+1}},
\end{eqnarray*}
where in the last step we have used the assumption that
$k+1-(l+1)\geq 0.$ So the identity holds
for $l+1.$ \qed

\end{document}